\documentclass{article}

 \usepackage{amsmath, amscd, amssymb, amsthm, bm}
 \usepackage{threeparttable}
 \usepackage{graphicx}
 \usepackage{color} 
\usepackage{pdfsync}  
\definecolor{dullmagenta}{rgb}{0.4,0,0.4}   
\definecolor{darkblue}{rgb}{0,0,0.4}
 \usepackage[colorlinks=true, pdfstartview=FitV, linkcolor=black, citecolor=black, urlcolor=black,plainpages=false,pdfpagelabels]{hyperref}  
\hypersetup{ pdfauthor = {Ameet S. Deshpande}, pdftitle = {New Fundamental solution for the time varying differential Riccati equation}, pdfkeywords = {Keywords for Searching}, pdfcreator = {pdfLaTeX with hyperref package}, pdfproducer = {pdfLaTeX}}
\usepackage{memhfixc}  
\usepackage{xspace}

\newcommand{\grad}{\nabla}

\newcommand{\cI}{{\mathcal I}}

\newcommand{\real} {\mathbb{R}}

\newcommand{\cIbarhatmk}{{\widehat{\overline \cI}^{\raise -0.3em \hbox{\scriptsize$m$}}_k}}
\newcommand{\cIbarhatmf}{{\widehat{\overline \cI}^{\raise -0.3em \hbox{\scriptsize$m$}}_f}}

\newcommand{\cIbartildemtk}{{\widetilde{\overline \cI}^{\raise -0.3em \hbox{\scriptsize$m,t$}}_k}}
\newcommand{\cIbartildemtf}{{\widetilde{\overline \cI}^{\raise -0.3em \hbox{\scriptsize$m,t$}}_f}}

\newcommand{\gamovertwo} {{\frac{\gamma^2}{2}}}

\newcommand{\gammovertwo} \gamovertwo
\newcommand{\gammuovertwo} \gamovertwo
\newcommand{\gammuepstovertwo} \gamovertwo
\newcommand{\gammubarepstovertwo} \gamovertwo

\newcommand{\Vttdel} {\widetilde{\widetilde{V}}^{\raise -0.3em \hbox{\scriptsize$\delta$}} }

\newcommand{\maxp} {{max-plus}\xspace}

\newcommand{\beasnum}{\begin{eqnarray}}
\newcommand{\eeasnum}{\end{eqnarray}}
\newcommand{\beas}{\begin{eqnarray*}}
\newcommand{\eeas}{\end{eqnarray*}}

\newcommand{\be}{\begin{equation}}
\newcommand{\ee}{\end{equation}}
\newcommand{\ba}{\begin{array}}
\newcommand{\ea}{\end{array}}
\newcommand{\er}[1]{\hbox{(\ref{#1})}}

\newtheorem{theorem}            {Theorem}[section]
\newtheorem{corollary}          [theorem]{Corollary}

\newtheorem{definition}         [theorem]{Definition}
\newtheorem{lemma}              [theorem]{Lemma}
\newtheorem{sideremark}         [theorem]{Remark}

\newtheorem{sideeg}           [theorem]{Example}
\newtheorem{sideconj}           [theorem]{Conjecture}

\newenvironment{remark}         {\begin{sideremark}\rm}{\end{sideremark}}

\newcommand{\til}[1]  {\tilde{#1}}

\newcommand{\bmat}{\left[\begin{matrix}}
\newcommand{\emat}{\end{matrix}\right]}

\title{\LARGE \bf
Max-Plus Representation for the Fundamental Solution of the Time-Varying Differential Riccati Equation
}

\author{Ameet Shridhar Deshpande
\thanks{Ameet Deshpande is with
Clipper windpower,
6305 Carpinteria Avenue,Carpinteria, CA 93013, USA
{\tt\small ameet.deshpande@gmail.com}}
\thanks{Research partially supported by NSF grant DMS-0307229 and
AFOSR grant FA9550-06-1-0238.}
}

\begin{document}

\maketitle
\thispagestyle{empty}

\pagestyle{empty}

\begin{abstract}
Using the tools of optimal control, semiconvex duality and \maxp algebra, this work derives a unifying representation of the solution for the matrix differential Riccati equation (DRE) with time-varying coefficients.
It is based upon a special case of the \maxp fundamental solution, first proposed in \cite{FlemMac}. Such fundamental solution can extend a special solution of certain bivariate DRE into the general solution, and can analytically solve the DRE starting from any initial condition.

This paper also shows that under a fixed duality kernel, the semiconvex dual of a DRE solution satisfies another dual DRE, whose coefficients satisfy the matrix compatibility conditions involving Hamiltonian and certain symplectic matrices.
For the time invariant DRE, this allows us to make dual DRE linear and thereby solve the primal DRE analytically. 
This paper also derives various kernel/duality relationships between the primal and time shifted dual DREs,
which leads to an array of DRE solutions. Time invariant analogue of one of these methods was first proposed in \cite{Funda}.
\end{abstract}
\section{Introduction}

The differential Riccati equation (DRE) plays a central role in estimation and optimal control.

An extensive study of algorithms for solving time-invariant and time-varying DREs was carried out by Kenney and Leipnik \cite{Leipnik}.
These include direct integration, the Chandrasekhar, Leipnik, Davison-Maki, modified Davison-Maki algorithms.
Later important developments include a Bernoulli substitution algorithm by Laub \cite{Laub82}, eigenvector decomposition techniques by
Oshman and Bar-Itzak \cite{oshman}, generalized partitioned solutions and integration free algorithms by Lainotis \cite{Lainiotis}, 
superposition laws developed by Sorine and Winternitz \cite{SorineWinternitz}, 
solutions by Rusnak \cite{Rusnak}, \cite{Rusnak98}. 
More recently, a fundamental solution based on \maxp algebra and semiconvex duality
was proposed by McEneaney \cite{Funda}.

%


The purpose of this paper is to present a new representation of the fundamental solution of the time-varying DRE.
The fundamental solution allows us to efficiently compute a general solution starting from different initial conditions.
This representation uses the \maxp techniques and is inspired from \cite{Funda}, but it extends the solution to the time-varying DRE and simplifies the treatment
by not using the semiconvex duality for the main result. In process, it derives the special case of the \maxp fundamental solution first proposed by Fleming and McEneaney in \cite{FlemMac}, for the linear-quadratic problem. It also shows that such fundamental solution is bivariate quadratic and describes the algorithm to compute the same.
It shows that evolution of a DRE under the \maxp fundamental solution is also a semiconvex dual transformation with a suitable kernel. Further it shows that the semiconvex dual transformation of a DRE, satisfies another DRE. It then derives the matching conditions between the coefficients and duality kernels relationships between primal and dual solutions at different times.

The DRE solution itself is similar in structure to the previous algorithms. Specifically, the fundamental solution computation requires integration of three ODEs similar to the forward formulas in \cite{Lainiotis} and $1$-representation addition formula in \cite{SorineWinternitz}. Still, the \maxp framework presented here is unifying and general.
E.g. partitioned formulas for the forward and backward time-varying DREs in \cite{Lainiotis}, time-invariant DRE solutions in \cite{Funda}, \cite{Leipnik2}, \cite{Rusnak} can be derived as special cases of a single framework. In addition, it is known that such algorithms work well for the stiff time-varying DREs and long time horizons without any computational difficulties, unlike the time-marching algorithms or the Davison-Maki algorithm.

\section{Optimal control problem}

We consider the matrix  \emph{differential Riccati equation} (DRE) of the form
\begin{align}\label{eq:ric:introdre}
\begin{split}
 -\dot{p}(t)&=A(t){'}p(t)+p(t)A(t)+C(t)+p(t){'}\Sigma(t) p(t)
 \end{split}
\end{align}
given the boundary condition $p(T)$ at time $T$. Here $t\in (-\infty,T]$ and $A(t)$ is square, $p(T),C(t),\Sigma(t)$ are square 
and symmetric $n\times n$ matrices and $\Sigma(t) = \sigma(t)\sigma(t){'} \succeq 0$ where $\sigma(t)$ is $n\times m$ matrix. 
It is well-known that above DRE arises in optimal control problem with linear dynamics 
\begin{align}
\dot{\xi}_{s}&= f_{s}(\xi_{s},u_{s})\doteq A(s){\xi}_{s}+\sigma(s){u}_{s}, \quad \xi_{t}=x\in \real^{n} \label{eq:ric:dyneq} \\
\intertext{and the following payoff function consisting of the integral and terminal payoffs,}
J_{t}^{T}(x,u) &\doteq \int_{t}^{T} l_{s}(\xi_{s},u_{s})\, ds + \phi(\xi_{T}), \text{\quad where} \label{eq:ric:payoff}\\
 l_{s}(\xi,\upsilon) &\doteq \frac{1}{2}\xi^{'}C(s)\xi - \frac{1}{2}|\upsilon|^{2} \text{\quad and} \label{eq:ric:runpayoff}\\
\phi(\xi) &\doteq \frac{1}{2}\xi^{'}p(T)\xi, \quad \text{for all $\xi\in\real^{n}, \upsilon\in\real^{m}$.
}\label{eq:ric:termpayoff}
\end{align}

\noindent Then, the optimal payoff or the value function is also quadratic given by,
\be \label{eq:ric:simpleV}
V(t,x) \doteq \sup_{u\in L_{2}(t,T)} J_{t}^{T}(x, u) = \frac{1}{2}x{'}p(t)x,
\ee 
 and $p(t)$ follows the DRE \er{eq:ric:introdre}.
In order to ensure the existence and the regularity of the value function and for the development to follow, we make following assumptions.
%
\vspace{0.05 in}
\be\label{eq:ric:assm1}
\begin{minipage}{0.85 \linewidth}
{
We assume that $A(t), C(t), \Sigma(t)$ are piecewise continuous, locally bounded functions of time $t$ and that $\Sigma(t)\doteq \sigma(t)\sigma(t)'\succeq 0$, $\forall t \in \left(\bar{T},T\right]$. We also assume that the underlying dynamic system \er{eq:ric:dyneq} is controllable. Since the DRE may exhibit finite time blowup, we assume that for $t\in \left(\bar{T},T\right]$ with $t \leq T$, there exists a solution of DRE \er{eq:ric:introdre} with the terminal condition $p(T)=P_{T}$. We denote this solution by $P_{t}$ for the ease of notation.
}
\end{minipage}
\ee
\vskip 0.5 cm

\noindent Now we shall obtain the fundamental solution for DRE \er{eq:ric:introdre} through the following generalization of the above optimal control problem. 
We assume the same dynamics as \er{eq:ric:dyneq}, and assume the following payoff function in which the integral payoff $l_{s}$ is as defined in \er{eq:ric:runpayoff} and 
the terminal payoff is parametrized by $z\in \real^{n}$.
\begin{align}
\begin{split} \label{eq:biric:payoffeq}
 {J_{t}^{T}}(x,u; z) &\doteq \int_{t}^{T}l_{s}(\xi_{s},u_{s})\, ds + \phi(\xi_{T}; z), \quad \text{where}\\
 \phi(\xi; z) &\doteq \phi^{z}(\xi) \doteq \frac{1}{2}\xi{'}P_{T}\xi+\xi{'}S_{T} z +\frac{1}{2}z{'}Q_{T} z,
\quad \forall \, \xi\in\real^{n}.
\end{split}
\end{align}

%

\noindent The optimal payoff or the value function is defined as
\be 
V(t,x; z)\doteq V_{t}(x; z)\doteq V_{t}^{z}(x) \doteq \sup_{u\in L_2(t,T)} {J_{t}^{T}}(x,u; z)
\label{eq:ric:valuedef}\ee
for all $x,z \in \real^n$ and $t \in \left(\bar{T},T\right]$.

Now we state two important lemmas regarding such a value function, which are proved in the appendix.

\begin{theorem}
\label{thm:biquadhjbexist}
Assume \er{eq:ric:assm1}, and assume that $P_{T}$, $Q_{T}$ are symmetric matrices and $S_{T}$ is invertible. 
Then for any $z\in \real^n$, the value function \er{eq:ric:valuedef} is given by.
\be\label{eq:ric:Vzdef}
V(t,x;z)=\frac{1}{2}x{'}P_t x+x{'}S_t z+\frac{1}{2}z{'}Q_t z
\ee
where $P_t$, $S_t$, $Q_t$ evolve as per
\beasnum
-\dot{P_t}&=&A(t){'}P_t+P_t A(t)+C(t)+P_t\Sigma(t) P_t \nonumber\\
-\dot{S_t}&=&(A(t)+\Sigma(t)P_t){'}S_t \label{eq:ric:biriccati}\\
-\dot{Q_t}&=&S_t{'}\Sigma(t)S_t, \nonumber
\eeasnum
and satisfy the boundary conditions $P_{T}$, $S_{T}$ and $Q_{T}$,  respectively, at time $t=T$. 
Further, the optimal control at a state $\til{\xi}_{s}$ at time $s$ is 
\be \label{eq:optcontrolz}
\til{u}_{s} = \sigma(s){'}(P_{s}\til{\xi}_{s}+S_{s}z),
\ee
and the corresponding optimal trajectory $\til{\xi}$, starting at $\til{\xi}_{t}=x$ and evolving as per the control \er{eq:optcontrolz} satisfies 
\be \label{eq:ric:opttrajend}
S_{t_{2}}{'}\tilde{\xi}_{t_{2}}+Q_{t_{2}}z=S_{t_{1}}{'}\tilde{\xi}_{t_{1}}+Q_{t_{1}}z,
\ee
for $\bar{T} < t_{1} < t_{2} \leq T$. Further, $Q_{t_{1}} - Q_{t_{2}} \succ 0$ and $S_{t}$ is invertible for $t\in (\bar{T},T]$.
\end{theorem}
\begin{proof}
Lemmas \ref{lemma:biquadhjbexist1}, \ref{lem:ric:verification}, \ref{lem:ric:opttrajend}, \ref{lem:ric:posdefQ} in the appendix, together prove the above result.
\end{proof}

\begin{remark}\label{rem:ric:optstartpt}
Since $S_{t_1}$ and $S_{t_2}$ are invertible, \er{eq:ric:opttrajend} suggests a one-one and onto relation between start and end of optimal trajectories, $\xi_{t_1}$ and $\xi_{t_2}$ for all $z$. Thus $\forall y \in \real^n$ there exists a $x=S_{t_2}^{-1}{'}\left( S_{t_1}{'}y+(Q_{t_1}-Q_{t_2})z\right)$ such that optimal trajectory $\tilde{x}_t$ starting at 
$\tilde{x}_{t_1}=x$, ends with $y$. Thus every $y\in \real^n$ is an optimal point for some initial condition. 
\end{remark}


\section{Max-Plus Fundamental Solution}

Given $t_{1}, t_{2}\in \real$, $t_{1}<t_{2}$, system trajectory starting at $\xi_{t_{1}}=x$ and a general terminal payoff function $\psi:\real^{n}\rightarrow \real$, let us define the operator, 

\be
\mathcal{S}_{t_{1}}^{t_{2}}[\psi](x) \doteq \sup_{u\in L_{2}(t_{1},t_{2})}\int_{t_{1}}^{t_{2}}l_{s}(\xi_{s},u_{s})\, ds+ \psi(\xi_{t_{2}})
\label{eq:ric:semigroup}\ee

We can restate \er{eq:ric:valuedef} and \er{eq:biric:payoffeq} using above operator. 
Noting that $V_{T}^{z}(x)=\phi^{z}(x)$, as defined in \er{eq:biric:payoffeq}, we have for all $t\in\left(\bar{T},T\right]$
$$V_{t}^{z}(x)=\mathcal{S}_{t}^{T}[\phi^{z}](x)=\mathcal{S}_{t}^{T}[V_{T}^{z}](x)$$
It is well known that operators $\mathcal{S}_{t_{1}}^{t_{2}}$ form a semigroup.  That is if $t_{1}\leq t \leq t_{2}\leq T$, then
$\mathcal{S}_{t_{1}}^{t_{2}}[\psi]=\mathcal{S}_{t_{1}}^{t}[\mathcal{S}_{t}^{t_{2}}[\psi]]$, which is the celebrated Dynamic programming principle for this problem. That is with $t_{2}=T$,
\begin{align}\label{eq:ric:DPP}
\begin{split}
V_{t_{1}}^{z}(x)&=\mathcal{S}_{t_{1}}^{T}[\phi^{z}](x)=\mathcal{S}_{t_{1}}^{t}\left[\mathcal{S}_{t}^{T}[\phi^{z}]\right](x)=\mathcal{S}_{t_{1}}^{t}\left[V_{t}^{z}\right](x)\\
&=\sup_{u\in L_{2}(t_{1},t)}\int_{t_{1}}^{t}l_{s}(x_{s},u_{s})\, ds+ V_{t}^{z}(\xi_{t})
\end{split}
\end{align}
If we define $a \oplus b \doteq \max(a,b)$ and $a \otimes b\doteq a+b$, then it is well known that $\left(\real \cup \{-\infty\}, \oplus,\otimes\right)$ forms a commutative semifield which is referred to as the \maxp algebra (see \cite{baccellietal},\cite{heltonjames}, \cite{litvmaslov} for a fuller discussion).
We can extend this algebra to functions so as to define the \maxp vector space. Let $[a\oplus b](x)=\max(a(x),b(x))$ and $a(x)\otimes k=a(x)+k$, where $a,b: \real^{n}\rightarrow \real$ and $k\in \real$.
Maslov \cite{maslov} proved that the above semigroup is linear in \maxp algebra. Thus using above notation
\begin{align*} 
\mathcal{S}_{t_{1}}^{T}[\psi_{1}\oplus \psi_{2}](x) &\doteq \mathcal{S}_{t_{1}}^{t_{2}}[\max(\psi_{1},\psi_{2})](x) \\
&=\max\left\{\mathcal{S}_{t_{1}}^{t_{2}}[\psi_{1}](x),\mathcal{S}_{t_{1}}^{t_{2}}[\psi_{2}](x)  \right\} \doteq\mathcal{S}_{t_{1}}^{T}[\psi_{1}](x)\oplus \mathcal{S}_{t_{1}}^{T}[\psi_{2}](x)\\
\intertext { and }
\mathcal{S}_{t_{1}}^{T}[k\otimes \psi_{1}](x) &\doteq \mathcal{S}_{t_{1}}^{t_{2}}[k+\psi_{1}](x)
=k+\mathcal{S}_{t_{1}}^{t_{2}}[\psi_{1}](x)
\doteq k\otimes \mathcal{S}_{t_{1}}^{T}[\psi_{1}](x).
\end{align*}

Now we shall define a {\it \maxp kernel} $I: \real^{n}\times \real^{n}\rightarrow \real$ derived earlier in \cite{FlemMac} and \cite{Flemingbook}.
Let $\bar{T}<t_{1}\leq t_{2} \leq T$ and $x,y \in \real^{n}$, and $\xi_{t}$ evolve with dynamics \er{eq:ric:dyneq}. Define
\begin{align} \label{eq:ric:Fundakernel}
I_{t_{1}}^{t_{2}}(x,y)&\doteq \left \{ 
\ba{ll}
\sup_{u\in\mathcal{U}_{t_1}^{t_2}(x,y)} \int_{t_{1}}^{t_{2}}l_{t}(\xi_{t},u_{t})\, dt  
& \mbox{if } \mathcal{U}_{t_1}^{t_2}(x,y)\neq \emptyset \\
-\infty & \mbox{otherwise}
\ea \right. \\
\intertext{where}
\mathcal{U}_{t_1}^{t_2}(x,y)&\doteq \left\{u\in L_{2}(t_{1},t_{2}) : \xi_{t_{1}}=x, \xi_{t_{2}}=y \right\} \label{eq:ric:2ptinpspace}
\end{align}
Note that $I_{t_{1}}^{t_{2}}=-\infty$ indicates that it is impossible to reach $y$ from $x$ in time interval $(t_{1},t_{2})$ using any possible control $u$.


Fleming and McEneaney \cite{FlemMac} proposed  above kernel, and showed that
\be \label{eq:ric:kernelop}
\mathcal{S}_{t_{1}}^{t_{2}}[\psi](x)=\sup_{y\in\real^{n}}\left( I_{t_{1}}^{t_{2}}(x,y)+ \psi(y) \right)\doteq \int_{\real^{n}}^{\oplus}I_{t_{1}}^{t_{2}}(x,y)\otimes \psi(y)\, dy
\ee
and since $I_{t_{1}}^{t_{2}}$ depends only on the dynamics $\dot{\xi_{s}}=f_{s}(\xi_{s},u_{s})$ and running payoff $l_{s}(\xi_{s},u_{s})$, it is independent of the terminal payoff $\psi$. Hence it can serve as a {\it Fundamental solution}, and obtain $\mathcal{S}_{t_{1}}^{t_{2}}[\psi](x)$ for any $\psi(x)$ by a kernel operation. 


\begin{remark}\label{rem:ric:Unonempty}
Also note that due to the controllability assumption \er{eq:ric:assm1}, for $t_{1}<t_{2}$, $\mathcal{U}_{t_1}^{t_2}(x,y)\neq \emptyset$ and for some $u \in \mathcal{U}_{t_1}^{t_2}(x,y)$,
 $I_{t_1}^{t_2}(x,y) \geq \int_{t_1}^{t_2} l_{s}(\xi_{s},u_{s})\, ds > -\infty$ for all $x,y\in\real^n$. Using \er{eq:ric:kernelop} and \er{eq:biric:payoffeq}, $\mathcal{S}_{t_{1}}^{t_{2}}[\phi^{z}](x) > -\infty$.
For $t_{1}=t_{2}$, $I_{t_{1}}^{t_{2}}(x,y)=-\infty$ for all $y\neq x$ and $I_{t_1}^{t_2}(x,x)=0$.
\end{remark}

\subsection{Computing the \maxp kernel}

\begin{theorem}
Assume \er{eq:ric:assm1}. Let $V$ and $I$ be as per \er{eq:ric:valuedef} and \er{eq:ric:Fundakernel}, respectively. 
Assume $x,y\in \real^n$ and $\bar{T}<t_1 < t_2 \leq T$. Thus using \er{eq:ric:Vzdef}, $V_{t_1}^z(\xi)=\frac{1}{2}\xi{'}P_{t_1}\xi+\xi{'}S_{t_1}z+\frac{1}{2}z'Q_{t_1}z$ and
$V_{t_2}^z(\xi)=\frac{1}{2}\xi{'}P_{t_2}\xi+\xi{'}S_{t_2}z+\frac{1}{2}z'Q_{t_2}z$. Then,
The \maxp kernel $I_{t_1}^{t_2}(x,y)$ exists and can be computed as
\be \label{eq:ric:computekernel}
\inf_{z\in\real^n}\left[ V_{t_1}^z(x)-V_{t_2}^z(y)\right]=V_{t_1}^{\hat{z}}(x)-V_{t_2}^{\hat{z}}(y)=I_{t_1}^{t_2}(x,y),
\ee
where $\hat{z}=(Q_{t_1}-Q_{t_2})^{-1}(S_{t_2}{'}y-S_{t_1}{'}x)$.
The \maxp kernel is also bivariate quadratic given by
\begin{align}\label{eq:ric:kernelform}
\begin{split}
I_{t_1}^{t_2}(x,y)&=\frac{1}{2}x'{I_{11}}_{t_1}^{t_2}x+x'{I_{12}}_{t_1}^{t_2}y+\frac{1}{2}y'{I_{22}}_{t_1}^{t_2}y,  \hskip 2em \text{where}\\
{I_{11}}_{t_1}^{t_2}&=P_{t_1}-S_{t_1}(Q_{t_1}-Q_{t_2})^{-1}S_{t_1}{'}\\
{I_{12}}_{t_1}^{t_2}&=S_{t_1}(Q_{t_1}-Q_{t_2})^{-1}S_{t_2}{'}\\
{I_{22}}_{t_1}^{t_2}&=-P_{t_2}-S_{t_2}(Q_{t_1}-Q_{t_2})^{-1}S_{t_2}{'}.
\end{split}
\end{align}
\end{theorem}
\begin{proof}
Let $\xi_{t_1}=x$. 
\begin{align}
&V_{t_1}^z(x)-V_{t_2}^z(y) =\mathcal{S}_{t_1}^{t_2}[V_{t_2}^z](x)-V_{t_2}^z(y) \nonumber\\
&\hskip 2em =\sup_{u\in L_2(t_1,t_2)}
\left\{\int_{t_1}^{t_2}l_{s}(\xi_{s},u_{s})\,ds+V_{t_2}^z(\xi_{t_2})-V_{t_2}^z(y)\right\}\nonumber\\
\intertext{substituting for $V_{t_2}^z$,}
&\hskip 2em =\sup_{u\in L_2(t_1,t_2)}
\left\{\int_{t_1}^{t_2}l_{s}(\xi_{s},u_{s})\,ds+\frac{1}{2}\xi_{t_2}'P_{t_2}\xi_{t_2}-\frac{1}{2}y'P_{t_2}y+ (\xi_{t_2}-y)'S_{t_2}z\right\}\nonumber\\
\intertext{Using \er{eq:ric:2ptinpspace}, $\mathcal{U}_{t_1}^{t_2}(x,y) \subset L_2(t_1,t_2)$ , and $\forall u \in \mathcal{U}_{t_1}^{t_2}(x,y)$, $\xi_{t_2}=y$.}
& \hskip 2em \geq \sup_{u\in \mathcal{U}_{t_1}^{t_2}(x,y)}\left\{\int_{t_1}^{t_2}l_{s}(\xi_{s},u_{s})\,ds+\frac{1}{2}y'P_{t_2}y-\frac{1}{2}y'P_{t_2}y+ (y-y)'S_{t_2}z\right\} \nonumber\\
& \hskip 2em = \sup_{u\in \mathcal{U}_{t_1}^{t_2}(x,y)}\int_{t_1}^{t_2}l_{s}(\xi_{s},u_{s})\,ds = I_{t_1}^{t_2}(x,y) \label{eq:ric:kernelform:b}
\end{align}
Taking infimum over all $z\in \real^n$,
\be\label{eq:ric:kernelform:c}
\inf_{z\in\real^n}\left[ V_{t_1}^z(x)-V_{t_2}^z(y)\right]\geq I_{t_1}^{t_2}(x,y)
\ee
Since $Q_{t_{1}}-Q_{t_{2}}\succ 0$ by theorem \ref{thm:biquadhjbexist}, define  $\hat{z}=(Q_{t_1}-Q_{t_2})^{-1}(S_{t_2}{'}y-S_{t_1}{'}x)$. Hence
$$S_{t_{2}}{'}y+Q_{t_{2}}\hat{z}=S_{t_{1}}{'}x+Q_{t_{1}}\hat{z}$$
hence using \er{eq:ric:opttrajend} the optimal trajectory $\tilde{\xi}_{s}$ starting from $\tilde{\xi}_{t_1}=x$ and with terminal payoff $V_{t_{2}}^{\hat{z}}(\cdot)$, ends at $\tilde{\xi}_{t_2}=y$. Let the corresponding optimal control be $\tilde{u}_{s}$.
Hence,
\begin{align}
V_{t_1}^{\hat{z}}(x)-V_{t_2}^{\hat{z}}(y)&=\left\{\sup_{u\in L_2(t_1,t_2)}\int_{t_1}^{t_2}l_{s}(\xi_{s},u_{s})\,ds+V_{t_2}^{\hat{z}}(\xi_{t_{2}})\right\}-V_{t_2}^{\hat{z}}(y) \nonumber\\
&=\left\{\int_{t_1}^{t_2}l_{s}(\til{\xi}_{s},\til{u}_{s})\,ds+V_{t_2}^{\hat{z}}(\til{\xi}_{t_{2}}) \right\} -V_{t_2}^{\hat{z}}(y) \nonumber\\
&=\int_{t_1}^{t_2}l_{s}(\til{\xi}_{s},\til{u}_{s})\,ds+V_{t_2}^{\hat{z}}(y)  -V_{t_2}^{\hat{z}}(y) \nonumber\\
\intertext{since $\tilde{u}\in\mathcal{U}_{t_1}^{t_2}(x,y)$}
&\leq \sup_{u\in \mathcal{U}_{t_1}^{t_2}(x,y)}\int_{t_1}^{t_2}l_{s}(\xi_{s},u_{s})\,ds
= I_{t_1}^{t_2}(x,y) \label{eq:ric:kernelform:e}
\end{align}
%

Inequalities \er{eq:ric:kernelform:c} and \er{eq:ric:kernelform:e} together give us \er{eq:ric:computekernel} and
%
substituting $\hat{z}$ in \er{eq:ric:computekernel} and expanding, we get \er{eq:ric:kernelform}.
\end{proof}

Now we shall prove a theorem which can allow us to combine \maxp kernels in time.
\begin{theorem}\label{thm:ric:combop}
Assuming \er{eq:ric:assm1}, let $\bar{T}<t_{1}<t_{2}<t_{3}\leq T$, then \maxp kernel $I_{t_{1}}^{t_{3}}$ can be computed from $I_{t_{1}}^{t_{2}}$ and $I_{t_{2}}^{t_{3}}$ as follows
\be\label{eq:ric:combop}
I_{t_{1}}^{t_{3}} (x,y)=\mathcal{S}_{t_{1}}^{t_{2}}[I_{t_{2}}^{t_{3}}(. , y)](x)=\sup_{z\in \real^{n}}\left\{ I_{t_{1}}^{t_{2}}(x,z)+I_{t_{2}}^{t_{3}}(z,y)\right\}
\ee
Thus $I_{t_1}^{t_3}(x,y)=\frac{1}{2}x{'}{I_{11}}_{t_1}^{t_3}x + x{'}{I_{12}}_{t_1}^{t_3}y + \frac{1}{2}y{'}{I_{22}}_{t_1}^{t_3}y$ where
\begin{align}\label{eq:ric:combop2}
\begin{split}
{I_{11}}_{t_1}^{t_3}&={I_{11}}_{t_1}^{t_2}-{I_{12}}_{t_1}^{t_2}
\left( {I_{22}}_{t_1}^{t_2}+{I_{11}}_{t_2}^{t_3}\right)^{-1}{{I_{12}}_{t_1}^{t_2}}^{'} \\
{I_{12}}_{t_1}^{t_3}&= -{I_{12}}_{t_1}^{t_2}
\left( {I_{22}}_{t_1}^{t_2}+{I_{11}}_{t_2}^{t_3} \right)^{-1}{I_{12}}_{t_2}^{t_3} \\
{I_{22}}_{t_1}^{t_3}&={I_{22}}_{t_2}^{t_3}-{{I_{12}}_{t_2}^{t_3}}^{T}
\left( {I_{22}}_{t_1}^{t_2}+{I_{11}}_{t_2}^{t_3}\right)^{-1} {I_{12}}_{t_2}^{t_3}
\end{split}
\end{align}
\end{theorem}
\begin{proof}
Note that by remark \ref{rem:ric:Unonempty}, $\mathcal{U}_{t_{a}}^{t_{b}}(x,y)\neq \emptyset$ for all $t_{a}<t_{b}$ and $x,y\in \real^{n}$.
\begin{align}
 I_{t_{1}}^{t_{3}} (x,y)&=\sup_{u\in\mathcal{U}_{t_{1}}^{t_{2}}(x,y)}\int_{t_{1}}^{t_{3}}l_{s}({\xi}_{s},u_{s})\,ds \nonumber\\
\intertext{ since $\mathcal{U}_{t_{1}}^{t_{3}}(x,y)=\bigcup_{z\in\real^n}\left(\mathcal{U}_{t_1}^{t_2}(x,z)\cap \mathcal{U}_{t_2}^{t_3}(z,y)  \right)$ } 
&= \sup_{z\in\real^n}\sup_{u\in U_{t_1}^{t_2}(x,z)\cap U_{t_2}^{t_3}(z,y)}\int_{t_{1}}^{t_{3}}l_{s}({\xi}_{s},u_{s})\,ds  \label{eq:ric:combop:a}
\end{align}
Now, we consider the following
\begin{align}
& \sup_{u\in U_{t_1}^{t_2}(x,z)\cap U_{t_2}^{t_3}(z,y)}\int_{t_{1}}^{t_{3}}l_{s}({\xi}_{s},u_{s})\,ds \nonumber\\
&\hskip 6em \leq \sup_{u\in U_{t_1}^{t_2}(x,z)\cap U_{t_2}^{t_3}(z,y)}\left\{\int_{t_{1}}^{t_{2}}l_{s}({\xi}_{s},u_{s})\,ds +\int_{t_{2}}^{t_{3}}l_{s}({\xi}_{s},u_{s})\,ds \right\} \nonumber\\
%
& \hskip 6em \leq \sup_{u\in U_{t_1}^{t_2}(x,z)}\left\{\int_{t_{1}}^{t_{2}}l_{t}({\xi}_{s},u_{s})\,ds\right\} +\sup_{u\in U_{t_2}^{t_3}(z,y)}\left\{\int_{t_{2}}^{t_{3}}l_{s}({\xi}_{s},u_{s})\,ds \right\}\nonumber\\
& \hskip 6em = I_{t_1}^{t_2}(x,z)+I_{t_2}^{t_3}(z,y) \label{eq:ric:combop:b}
\end{align}
Now, since $I_{t_1}^{t_2}(x,z) > -\infty$, $\forall \epsilon > 0$, $\exists$ $\bar{u} \in \mathcal{U}_{t_1}^{t_2}(x,z)$ and trajectory $\bar{\xi}_{s}$ with $\bar{\xi}_{t_1}=x$ such that 
\be\int_{t_1}^{t_2}l_{s}(\bar{\xi}_{s},\bar{u}_{s})\,ds +\epsilon \geq I_{t_1}^{t_2}(x,z)\label{eq:ric:combop:c}\ee 
Similarly $\exists$ $\tilde{u} \in \mathcal{U}_{t_2}^{t_3}(z,y)$ and trajectory $\tilde{\xi}$ with $\tilde{\xi}_{t_2}=z$ such that 
\be\int_{t_2}^{t_3}l_{s}(\tilde{\xi}_{s},\tilde{u}_{s})\,ds +\epsilon \geq I_{t_2}^{t_3}(z,y)\label{eq:ric:combop:d}\ee
Now we can create augmented control $\hat{u}$ such that $\hat{u}_{s}=\bar{u}_{s}$ for $s \in [t_1,t_2)$ and $\hat{u}_{s}=\tilde{u}_{s}$ for $t \in [t_2,t_3]$, and extend it arbitrarily beyond. Note that if $\hat{\xi}_{s}$ is corresponding trajectory, then starting with $\hat{\xi}_{t_1}=x$, $\hat{\xi}_{s}=\bar{\xi}_{s}$ for $t\in[t_1,t_2]$. Hence $\hat{\xi}_{t_2}=z$ and $\hat{\xi}_{s}=\tilde{\xi}_{s}$ for $s\in[t_2,t_3]$, hence $\hat{\xi}_{t_3}=y$. Hence $\hat{u}\in \mathcal{U}_{t_1}^{t_2}(x,z)\cap \mathcal{U}_{t_2}^{t_3}(z,y)$. Moreover using
\er{eq:ric:combop:c} and \er{eq:ric:combop:d},
\begin{align}
\sup_{u\in U_{t_1}^{t_2}(x,z)\cap U_{t_2}^{t_3}(z,y)}\int_{t_{1}}^{t_{3}}l_{s}({\xi}_{s},u_{s})\,ds 
&\geq \int_{t_{1}}^{t_{2}}l_{s}(\hat{\xi}_{s},\hat{u}_{s})\,ds +\int_{t_{2}}^{t_{3}}l_{s}(\hat{\xi}_{s},\hat{u}_{s})\,ds \nonumber\\
& = \int_{t_1}^{t_2}l_{s}(\bar{\xi}_{s},\bar{u}_{s})\,ds + \int_{t_2}^{t_3}l_{s}(\tilde{\xi}_{s},\tilde{u}_{s})\,ds \nonumber\\
& = I_{t_1}^{t_2}(x,z)+I_{t_2}^{t_3}(z,y) -2\epsilon \label{eq:ric:combop:e}
\end{align}
Since $\epsilon$ is arbitrary, from \er{eq:ric:combop:b}  and \er{eq:ric:combop:e}, we have 
$$\sup_{u\in U_{t_1}^{t_2}(x,z)\cap U_{t_2}^{t_3}(z,y)}\int_{t_{1}}^{t_{3}}l_{s}({\xi}_{s},u_{s})\,ds = I_{t_1}^{t_2}(x,z)+I_{t_2}^{t_3}(z,y) $$
which with \er{eq:ric:combop:a} proves \er{eq:ric:combop}.
Now, using \er{eq:ric:kernelform} and since $(Q_{t_1}-Q_{t_2})\succ 0$ and $(Q_{t_2}-Q_{t_3})\succ 0$,
\begin{align*}
&{I_{22}}_{t_1}^{t_2}+{I_{11}}_{t_2}^{t_3}\\
&\hskip 4em =\left(-P_{t_2}-S_{t_2}(Q_{t_1}-Q_{t_2})^{-1}S_{t_2}{'}\right)+ \left(P_{t_2}-S_{t_2}(Q_{t_2}-Q_{t_3})^{-1}S_{t_2}{'}\right)\\
&\hskip 4em =-S_{t_2}\left((Q_{t_1}-Q_{t_2})^{-1}+(Q_{t_2}-Q_{t_3})^{-1}\right)S_{t_2}{'}\\
&\hskip 4em \prec 0
\end{align*}
Thus $\left\{ I_{t_{1}}^{t_{2}}(x,z)+I_{t_{2}}^{t_{3}}(z,y)\right\}$ is concave in $z$. Thus supremum in \er{eq:ric:combop} exists, and we get \er{eq:ric:combop2} by algebraic computation of the local maxima.
\end{proof}
\begin{remark}
Note that $I_{t}^{t_3}(x,z)$ has the same bivariate form as $V_{t}^z$ given by \er{eq:ric:Vzdef}, and both $I_{t}^{t_3}$ and $V_{t}$ evolve in the time interval $(t_1,t_2)$ according to the semigroup $\mathcal{S}_{t_1}^{t_2}$ as per \er{eq:ric:combop}. Hence the parameters satisfy DREs similar to the \er{eq:ric:biriccati}.
\begin{align}\label{eq:ric:biriccatiI}
\begin{split}
-\frac{d}{dt}{{I}_{11}}_{t}^{t_3}&=A(t){'}{I_{11}}_{t}^{t_3}+{I_{11}}_{t}^{t_3} A(t)+C(t)+{I_{11}}_{t}^{t_3}\Sigma(t) {I_{11}}_{t}^{t_3} \\
-\frac{d}{dt}{I_{12}}_{t}^{t_3}&=(A(t)+\Sigma(t){I_{11}}_{t}^{t_3}){'}{I_{12}}_{t}^{t_3} \\
-\frac{d}{dt}{I_{22}}_{t}^{t_3}&={{I_{12}}_{t}^{t_3}}{'}\Sigma(t){I_{12}}_{t}^{t_3} 
\end{split}
\end{align}

\end{remark}

Next we discuss a way to compute the \maxp kernel from the transition matrix of the Hamiltonian.
Role of such transition matrix in DRE solutions is well-known, e.g. Davison-Maki method \cite{DavisonMaki}.
\begin{corollary}
Assume \er{eq:ric:assm1}. Define the $2n\times 2n$ Hamiltonian matrix as
\begin{align} \label{eq:defHamilt}
\mathcal{H}_{t} \doteq \begin{bmatrix} A(t) & \Sigma(t)\\ -C(t) & -A(t){'} \end{bmatrix} 
\end{align}
Let the state transition matrix $\Phi$, associated with the linear time varying system $\dot{\xi}_{t} = -\mathcal{H}_{t} \xi_{t}$, consist of four $n\times n$ sub-matrices
\begin{align} \label{eq:defPhi}
\Phi(t_{2},t_{1}) \doteq 
\begin{bmatrix} \Phi_{11} & \Phi_{12} \\ \Phi_{21} & \Phi_{22} \end{bmatrix}.
\end{align}
Then parameters of the \maxp kernel $I_{t_{1}}^{t_{2}}(x,y)$ per \er{eq:ric:kernelform}, are given by
\begin{align}\label{eq:ric:IfromPhi}
\begin{split}
{I_{11}}_{t_1}^{t_2}&=\Phi_{21}{\Phi_{11}}^{-1}
+{\Phi_{11}}^{-1}{'}\Phi_{12}^{-1}\\
{I_{12}}_{t_1}^{t_2}&=-{\Phi_{11}}^{-1}{'}\Phi_{12}^{-1}{\Phi_{11}}\\
{I_{22}}_{t_1}^{t_2}&=\Phi_{12}^{-1}{\Phi_{11}}
\end{split}
\end{align}
\end{corollary}
\begin{proof}
It can be verified that differential equations \er{eq:ric:biriccati} are equivalent to the 
following single matrix differential equation involving Hamiltonian matrix $\mathcal{H}$ and a symplectic matrix $\mathcal{K}$.
\begin{align}
-\dot{\mathcal{K}_{t}}=\mathcal{K}_{t}\mathcal{H}_{t} {\quad \text{where }}
\mathcal{K}_{t} \doteq \begin{bmatrix} S_{t}^{-1}P_{t} & -S_{t}^{-1}\\ S_{t}{'}-Q_{t}S_{t}^{-1}P_{t} & Q_{t}S_{t}^{-1} \end{bmatrix}.
\label{eq:ric:hamilmatch1}
\end{align}
This being a linear ODE, the solution is given by 
$$\mathcal{K}_{t_{2}}=\mathcal{K}_{t_{1}} \Phi(t_{2},t_{1})$$
\begin{align*}
\begin{bmatrix} S_{t_{2}}^{-1}P_{t_{2}} & -S_{t_{2}}^{-1}\\ S_{t_{2}}{'}-Q_{t_{2}}S_{t_{2}}^{-1}P_{t_{2}} & Q_{t_{2}}S_{t_{2}}^{-1} \end{bmatrix}
=\begin{bmatrix} S_{t_{1}}^{-1}P_{t_{1}} & -S_{t_{1}}^{-1}\\ S_{t_{1}}{'}-Q_{t_{1}}S_{t_{1}}^{-1}P_{t_{1}} & Q_{t_{1}}S_{t_{1}}^{-1} \end{bmatrix}
\begin{bmatrix} \Phi_{11} & \Phi_{12}\\ \Phi_{21} & \Phi_{22} \end{bmatrix} 
\end{align*}

Matching terms, we get following set of equations
\begin{align} \label{eq:ric:SQPfrmHamil}
\begin{split}
S_{t_{1}}&=(\Phi_{11}+\Phi_{12}P_{t_{2}})^{-1}{'}S_{t_{2}}\\
Q_{t_{1}}&=Q_{t_{2}}-S_{t_{2}}{'}(\Phi_{11}+\Phi_{12}P_{t_{2}})^{-1}\Phi_{12}S_{t_{2}}\\
P_{t_{1}}&=\left( \Phi_{21}+\Phi_{22}P_{t_{2}}\right)
\left( \Phi_{11}+\Phi_{12}P_{t_{2}}\right)^{-1}
\end{split}
\end{align}

Substituting \er{eq:ric:SQPfrmHamil} in \er{eq:ric:kernelform} gives us
\begin{align*}
{I_{11}}_{t_1}^{t_2}&=(\Phi_{21}+\Phi_{22}P_{t_2})(\Phi_{11}+\Phi_{12}P_{t_2})^{-1}
+(\Phi_{11}+\Phi_{12}P_{t_2})^{-1}{'}\Phi_{12}^{-1}\\
{I_{12}}_{t_1}^{t_2}&=-(\Phi_{11}+\Phi_{12}P_{t_2})^{-1}{'}\Phi_{12}^{-1}(\Phi_{11}+\Phi_{12}P_{t_2})\\
{I_{22}}_{t_1}^{t_2}&=\Phi_{12}^{-1}\Phi_{11}
\end{align*}
But since $(I_{11},I_{12},I_{22})$ depend only on $(A(t),C(t),\Sigma(t),t_1,t_2)$ and are independent of starting $(P,S,Q)$.
Thus above equations hold true for any $P_{t_{2}}$. Specifically, we can take $P_{t_2}=0$ to get \er{eq:ric:IfromPhi}.

\end{proof}

\section{Solving the DRE}

Now we shall study how the \maxp kernel can be used to solve the DRE \er{eq:ric:introdre}.
\begin{corollary}\label{corr:kernel_to_dresol}
Assume \er{eq:ric:assm1}. Let $\bar{T}<t_{1}<t_{2}\leq T$. 
Also assume that a solution to the DRE \er{eq:ric:introdre} with a possibly different terminal condition, $p(T)=\bar{p}_{T}$, exists for all $s\in[t_{1},T]$.

If we denote such solution by $\bar{p}$, then the solution $\bar{p}_{t_{1}}$, can be computed from $\bar{p}_{t_{2}}$ as per
\be \label{eq:kernel_to_dresol:a}
\bar{p}_{t_1} = I_{11}-I_{12}(\bar{p}_{t_2}+I_{22})^{-1}I_{12}{'},
\ee
where $I_{11}$, $I_{12}$ and $I_{22}$ are the parameters of the kernel $I_{t_{1}}^{t_{2}}$.

Moreover if $(P_{t},S_{t},Q_{t})$ are the particular solutions of coupled DREs \er{eq:ric:biriccati} satisfying suitable boundary conditions from theorem \ref{thm:biquadhjbexist}, then following is an alternate method of propagating $\bar{p}$
\be\label{eq:ric:maxpluspropag2}
\bar{p}_{t_{1}}=P_{t_{1}}-S_{t_{1}}\left( Q_{t_{1}}-Q_{t_{2}}-S_{t_{2}}{'}(\bar{p}_{t_{2}}-P_{t_{2}})^{-1}S_{t_{2}}\right)^{-1}S_{t_{1}}{'}.
\ee
\end{corollary}
\begin{proof}
Using theorem \ref{eq:ric:computekernel}, the \maxp kernel $I_{t_{1}}^{t_{2}}$ exists.
With the terminal payoff $\phi(\xi) \doteq \frac{1}{2}\xi^{'}\bar{p}_{T}\xi$ and integral payoff
defined in \er{eq:ric:runpayoff}, let $V$ and $\mathcal{S}$ be the corresponding value function and the semigroup defined as per \er{eq:ric:simpleV} and \er{eq:ric:semigroup} respectively. Then using \er{eq:ric:kernelop} and \er{eq:ric:DPP},
\begin{align*}
V_{t_{1}}(x) = \mathcal{S}_{t_{1}}^{T} [\phi](x) = \mathcal{S}_{t_{1}}^{t_{2}} [V_{t_{2}}](x)
= \sup_{y\in\real^{n}}\left( I_{t_{1}}^{t_{2}}(x,y)+ V_{t_{2}}(y) \right).
\end{align*}
Using $V_{t_{1}}(x)=\frac{1}{2}x{'}\bar{p}_{t_{1}}x$ and $V_{t_{2}}(x)=\frac{1}{2}x{'}\bar{p}_{t_{2}}x$,
substituting the parameters of $I_{t_{1}}^{t_{2}}$ from \er{eq:ric:kernelform}, we get
\begin{align}\label{eq:ric:premaxpluspropag}
\frac{1}{2}x{'}\bar{p}_{t_1}x &= \sup_y \left\{ \frac{1}{2}x'{I_{11}}x+x'{I_{12}}y+\frac{1}{2}y'({I_{22}}+\bar{p}_{t_2})y \right\}.
\end{align}
Matching the terms gives us \er{eq:kernel_to_dresol:a}.

In order to prove \er{eq:ric:maxpluspropag2}, substitute the kernel parameters from \er{eq:ric:kernelform} in \er{eq:kernel_to_dresol:a}. With some manipulation and Woodbury's matrix identity \cite{woodbury}, we get
\begin{align*}
\left(S_{t_1}{'}(\bar{p}_{t_1}-P_{t_1})^{-1}S_{t_1}\right)^{-1}
&=\left(-(Q_{t_1}-Q_{t_2})+S_{t_2}{'}(\bar{p}_{t_2}-P_{t_2})^{-1}S_{t_2}\right)^{-1} \\
\intertext{and therefore}
S_{t_1}{'}(\bar{p}_{t_1}-P_{t_1})^{-1}S_{t_1}+Q_{t_1}
&=S_{t_2}{'}(\bar{p}_{t_2}-P_{t_2})^{-1}S_{t_2}+Q_{t_2},
\end{align*}
which after rearrangement gives us \er{eq:ric:maxpluspropag2}.
\end{proof}
\begin{remark}
Note that we assumed that the propagation $\frac{1}{2}x{'}\bar{p}_{t_1}x=\mathcal{S}_{t_1}^{t_2}[V_{t_2}](x)$ exists, and derived  \er{eq:kernel_to_dresol:a}. This is also equivalent to ${I_{22}}+\bar{p}_{t_2}\prec 0$, so that the supremum in \er{eq:ric:premaxpluspropag} exists. Thus the set $\left\{ \bar{p}_{t_{2}} | \quad{I_{22}}_{t_1}^{t_2}+\bar{p}_{t_2}\prec 0 \right\}$ characterizes all the terminal conditions for which the solution exists (is norm bounded) for $t\in[t_{1},t_{2}]$.
Also note that the minimum time $\bar{T}$ for which solution to DRE exists, depends on initial condition. 
There are generalization to the DRE solution which allow solutions containing singularity e.g. \cite{SorineWinternitz}. Above propagation formula \er{eq:kernel_to_dresol:a} contains a pole, hence it is conjectured that it may allow computation of generalized DRE solutions which can step through singularity.
\end{remark}

\begin{remark}
Special cases of the alternate propagation formula \er{eq:ric:maxpluspropag2} yield both the forward and backward generalized partitioned formulae to solve DREs proposed in \cite{Lainiotis}. For small time propagation, when 
$t_{2} \rightarrow t_{1}$, $(Q_{t_1}-Q_{t_2})^{-1}$ and parameters of kernel $I_{t_1}^{t_2}$ defined in \er{eq:ric:kernelform} become more and more singular causing numerical inaccuracies in propagation. In such case, above alternate formula is useful, as it does not involve taking inverses of ill-conditioned matrices.
\end{remark}

\subsection{Algorithm}
Thus following is the final algorithm to obtain the fundamental solution, and to convert a particular solution of the  \er{eq:ric:biriccati} into the \maxp fundamental solution and a general solution of \er{eq:ric:introdre}. It gives us a closed form solution to the DRE \er{eq:ric:introdre} using  \maxp kernel $I_{t_1}^{t_2}$ \er{eq:ric:kernelop}. We shall reiterate the formulae derived earlier to make the section self-contained.

\begin{enumerate}

\item Compute the parameter triplet $(I_{11},I_{12},I_{22})_{t_1}^{t_2}$ of the {\it \maxp kernel} $I_{t_1}^{t_2}(x,y)$ defined in \er{eq:ric:kernelform} using any of the following methods.

\begin{enumerate}
\item Evolution of the bivariate payoff function using $3$ ODEs: 

Choose the parameters $(P_{t_{2}},S_{t_{2}},Q_{t_{2}})$ of the terminal bivariate payoff $V_{t_{2}}^{z}(x)=\frac{1}{2}x{'}P_{t_{2}}x+x{'}S_{t_{2}}z+\frac{1}{2}z{'}Q_{t_{2}}z$, such that $P_{t_{2}}$, $Q_{t_{2}}$ are $n\times n$ symmetric matrices, and $S_{t_{2}}$ is $n\times n$ invertible matrix.

Propagate $(P,S,Q)$ backwards in time according to \er{eq:ric:biriccati} till time $t_1 < t_2$. That is
\begin{align*}
-\dot{P_t}&=A(t){'}P_t+P_t A(t)+C(t)+P_t\Sigma(t) P_t\\
-\dot{S_t}&=(A(t)+\Sigma(t)P_t){'}S_t\\
-\dot{Q_t}&=S_t{'}\Sigma(t)S_t
\end{align*} 

Compute the parameters of the \maxp kernel,
\begin{align*}
{I_{11}}_{t_{1}}^{t_{2}}&=P_{t_1}-S_{t_1}(Q_{t_1}-Q_{t_2})^{-1}S_{t_1}{'}\\
{I_{12}}_{t_{1}}^{t_{2}}&=S_{t_1}(Q_{t_1}-Q_{t_2})^{-1}S_{t_2}{'}\\
{I_{22}}_{t_{1}}^{t_{2}}&=-P_{t_2}-S_{t_2}(Q_{t_1}-Q_{t_2})^{-1}S_{t_2}{'}.
\end{align*}

\item The state transition matrix of the Hamiltonian matrix:
Let the Hamiltonian matrix $\mathcal{H}_{t}$ be defined as per \er{eq:defHamilt}.
Let $\Phi_{11}$, $\Phi_{12}$, $\Phi_{22}$ be the submatrices of the $2n\times 2n$ state transition matrix $\Phi(t_{2},t_{1})$ be as per \er{eq:defPhi}. Then the paramteres of the \maxp kernel are computed as per \er{eq:ric:IfromPhi}. That is
\begin{align*}
\begin{split}
{I_{11}}_{t_1}^{t_2}&=\Phi_{21}{\Phi_{11}}^{-1}
+{\Phi_{11}}^{-1}{'}\Phi_{12}^{-1}\\
{I_{12}}_{t_1}^{t_2}&=-{\Phi_{11}}^{-1}{'}\Phi_{12}^{-1}{\Phi_{11}}\\
{I_{22}}_{t_1}^{t_2}&=\Phi_{12}^{-1}{\Phi_{11}}.
\end{split}
\end{align*}

\end{enumerate}
\item Given any terminal condition $\bar{p}_{t_{2}}$ at time $t_{2}$, if the solution exists, the DRE \er{eq:ric:introdre} can be solved at time $t_{1}$ to get $\bar{p}_{t_{1}}$ as per \er{eq:kernel_to_dresol:a}. That is
\begin{align*}
\bar{p}_{t_1}&=I_{11}-I_{12}(\bar{p}_{t_2}+I_{22})^{-1}I_{12}{'}.
\end{align*}
\item For small time propagation, as the parameters of kernel $I_{t_1}^{t_2}$ become ill-conditioned, the alternate formula
\er{eq:ric:maxpluspropag2} can be used. Instead of using the \maxp kernel, it uses the particular solution\er{eq:ric:biriccati} at times $t_{1}$ and $t_{2}$ directly
and avoids the inverses of near-singular matrices. That is
$$
\bar{p}_{t_{1}}=P_{t_{1}}-S_{t_{1}}\left( Q_{t_{1}}-Q_{t_{2}}-S_{t_{2}}{'}(\bar{p}_{t_{2}}-P_{t_{2}})^{-1}S_{t_{2}}\right)^{-1}S_{t_{1}}{'}
$$
This formula does not blow up for a small time step propagation, and yields an accurate propagation.
\end{enumerate}

\section{Semiconvex dual extensions}
Theory of semiconvex duality offers a natural language to express the \maxp fundamental solution derived before,
as the \maxp kernel operation can be expressed as a semiconvex dual operation with an appropriate kernel.
It has been used in \cite{Funda} to derive the semiconvex dual \maxp fundamental solution for time-invarient DREs.
Here we study matching conditions and various kernel relationships between primal DRE \er{eq:ric:introdre} and its dual,
which add to the toolchest to solve DREs.
In process, we also derive time varying extension of the fundamental solution proposed in \cite{Funda}.

\subsection{Semiconvex duality}
\begin{definition}
A function $\mathcal{P}(x):\real^n\rightarrow \real^- \doteq \real \cup \{-\infty\}$ is defined to be \emph{uniformly semiconvex} with (symmetric) matrix constant $K$ if $\mathcal{P}(x)+\frac{1}{2}x{'}Kx$ is strictly convex over $\real^n$. We denote this space by ${\mathcal{S}}^K$.
\end{definition}

Semiconvex duality is parametrized by a bivariate quadratic kernel 
\be\label{eq:ric:dualitykernel}
\phi(x,z)=\frac{1}{2}x{'}Px+x{'}Sz+\frac{1}{2}z{'}Qz
\ee
where $P$ and $Q$ are symmetric matrices. We use this kernel to define semiconvex duality.
\begin{lemma}\label{thm:ric:semicon}
Let $\mathcal{P} \in {\mathcal{S}}^{-P}$, $S$ is invertible and $\phi(x,z)$
defined as above. Then $\forall z \in \real^n$ we can define the duality operator $\mathcal{D}_{\phi}$ and the dual $\mathcal{Q}(z)$ of primal $\mathcal{P}(x)$ as follows.

\be 
\mathcal{D}_{\phi}[\mathcal{P}](z) \doteq  \inf_{x}[\mathcal{P}(x)-\phi(x,z)] \doteq \mathcal{Q}(z) \label{eq:ric:semicon}\ee 
from the dual $\mathcal{Q}(z)$, primal can be 
recovered again using the inverse duality operator $\mathcal{D}_{\phi}^{-1}$ defined below.
\be 
\mathcal{D}_{\phi}^{-1}[\mathcal{Q}](x)  \doteq \sup_{z}[\phi(x,z)+\mathcal{Q}(z)] \doteq \mathcal{P}(x). \label{eq:ric:semiconinv}
\ee
$\phi(x,z)$ is called the kernel of duality. Thus $\mathcal{D}_{\phi}^{-1}\mathcal{D}_{\phi}[\mathcal{P}](x)=\mathcal{P}(x)$.
\end{lemma}
\begin{proof}
Proved in appendix \ref{proof:thm:ric:semicon}.
\end{proof}
Using a very similar approach we can derive the following corollary. Here we state it without proof.
\begin{corollary}\label{corr:dualinv}
If $\mathcal{Q}(z)+z{'}Qz$ is concave over $z\in \real^{n}$, then 
\be\label{eq:ric:dualinv}
\mathcal{D}_{\phi}\mathcal{D}_{\phi}^{-1}\mathcal{Q}(z)=\mathcal{Q}(z)
\ee
\end{corollary}

If we choose $\mathcal{P}(x)=\frac{1}{2}x{'}px$, $\phi(x,z)=\frac{1}{2}x{'}Px+x{'}Sz
+\frac{1}{2}z{'}Q z$ and assume $p>P$ and $S$ is nonsingular, then $\mathcal{P}(x)\in \mathcal{S}^{-p}$. Hence by substitution in \er{eq:ric:semicon}, we get $\mathcal{Q}(z)=\frac{1}{2}z{'}qz$, where
\be q=-S{'}\left(p-P\right)^{-1}S-Q\label{eq:ric:semicon2}\ee
We can also derive the following inverse relation 
\be p=-S\left(q+Q\right)^{-1}S{'}+P\label{eq:ric:semiconinv2} \ee
\subsection{Dual DRE and Compatibility Conditions}
Now we consider a special duality problem involving only quadratic primal and dual functions.
Specifically, we assume the following.
\be\label{eq:ric:assm2}
\begin{minipage}{0.85 \linewidth}
{
Assume $\bar{T}<t_{1}<t_{2}\leq T$. Let $\bar{p}_{t}$ be the solution of \er{eq:ric:introdre} with the boundary condition $p_{t}=\bar{p}_{t_{2}}$.
Assume that such solution exists for $t_{1}\leq t \leq t_{2}$ and let $\mathcal{P}_t(x)\doteq \frac{1}{2}x{'}\bar{p}_t(x)$.
Recall the bivariate quadratic function in \er{eq:ric:Vzdef}. That is $V(t,x;z) = \frac{1}{2}x{'}P_{t}x+x{'}S_{t}z+\frac{1}{2}z{'}Q_{t}z$, with parameters $(P_{t},S_{t},Q_{t})$ evolving as per \er{eq:ric:biriccati}.
Also assume that $S_{t_{2}}$ is invertible.
\vskip 0.5em
Let the kernel of duality be $V_{t_{2}}(x,z)\doteq V(t_{2},x;z)$.
This kernel operates on $\mathcal{P}_{t}(x)$ to yield the dual $\mathcal{D}_{V_{t_{2}}}[\mathcal{P}_{t}](z)=\mathcal{Q}_{t}(z)=\frac{1}{2}z{'}\bar{q}_{t}(z)$ for all $t\in [t_{1},t_{2}]$.
Further assume that $S_{t_{2}}$ invertible, $\bar{p}_{t} \succ P_{t_{2}}$ for all $t \in [t_{1},t_{2}]$
and $P_{t_{1}} \succ P_{t_{2}}$.
}
\end{minipage}
\ee
\begin{theorem}
Assume \er{eq:ric:assm1} and \er{eq:ric:assm2}. Then $\bar{q}_{t}$ which parametrizes the semiconvex dual 
$\mathcal{Q}_{t}(z)=\mathcal{D}_{V_{t_{2}}}[\mathcal{P}_{t}](z)=\frac{1}{2}z{'}\bar{q}_{t}z$,
follows the dual DRE below
\be\label{eq:ric:ricd}
 -\dot{\bar{q}}_{t}=\bar{A}(t){'}\bar{q}_{t}+\bar{q}_{t}\bar{A}(t)+\bar{C}(t)_{}+\bar{q}_{t}{'}\bar{\Sigma}(t) \bar{q}_{t}.
\ee
The coefficients $(\bar{A},\bar{C},\bar{\Sigma})_{t}$ of the dual DRE satisfy the following matrix compatibility conditions.
Note that the time dependence of $(P, S, Q)_{t_{2}}$ is suppressed for brevity.
\begin{align}
\mathcal{K}_{t_{2}}\mathcal{H}_{t}&=\bar{\mathcal{H}}_{t}\mathcal{K}_{t_{2}}
\label{eq:ric:hamilmatch} \\
\intertext{where}
\mathcal{K}_{t_{2}}=\begin{bmatrix} S^{-1}P & -S^{-1}\\ S{'}-QS^{-1}P & QS^{-1} \end{bmatrix} ,
\mathcal{H}_{t}&=\begin{bmatrix} A(t)  & \Sigma(t) \\ -C(t) & -A(t){'} \end{bmatrix} ,
\bar{\mathcal{H}}_{t}=\begin{bmatrix} \bar{A}(t)  & \bar{\Sigma}(t) \\ -\bar{C}(t) & -\bar{A}(t){'} \end{bmatrix} \nonumber
\end{align}
 
\end{theorem}
\begin{proof}
Note that our assumptions allow us to use theorem \ref{thm:ric:semicon}.
Specifically, using \er{eq:ric:semicon2} and \er{eq:ric:semiconinv2}, we get
\be \bar{q}_{t}=-S_{t_{2}}{'}\left(\bar{p}_{t}-P_{t_{2}}\right)^{-1}S_{t_{2}}-Q_{t_{2}} \label{eq:ric:dualrica1}\ee
\be \bar{p}_{t}=-S_{t_{2}}\left(\bar{q}_{t}+Q_{t_{2}}\right)^{-1}S_{t_{2}}{'}+P_{t_{2}} \label{eq:ric:dualrica2}\ee
Differentiating both sides of \er{eq:ric:dualrica1},
\be \dot{\bar{q}}_{t}=S_{t_{2}}{'}\left(\bar{p}_{t}-P_{t_{2}}\right)^{-1}\dot{\bar{p}}_{t}\left(\bar{p}_{t}-P_{t_{2}}\right)^{-1}S_{t_{2}}\label{eq:ric:dualricb} \ee
Substituting for $\dot{\bar{p}}_{t}$ from \er{eq:ric:introdre}, $\bar{p}_{t}$ from \er{eq:ric:dualrica2} in \er{eq:ric:dualricb} and after simplification using \er{eq:ric:biriccati}, we get
 %
\begin{align*} 
-\dot{\bar{q}}_{t}&=\bar{q}_{t}{S_{t_{2}}}^{-1}\hat{P}_{t_{2}}{S_{t_{2}}^{-1}}^{'} \bar{q}_{t} 
+ \bar{q}_{t}{S_{t_{2}}}^{-1}( \hat{P}_{t_{2}}{S_{t_{2}}^{-1}}^{'} Q_{t_{2}}- \hat{S}_{t_{2}}) \\
&\hskip 0.4em + ( \hat{P}_{t_{2}}{S_{t_{2}}^{-1}}^{'} Q_{t_{2}}- \hat{S}_{t_{2}}){'}{{S_{t_{2}}}^{-1}}{'}\bar{q}_{t} 
+Q_{} {S_{t_{2}}}^{-1}\hat{P}_{t_{2}}{{S_{t_{2}}}^{-1}}{'}Q_{t_{2}} \\
&\hskip 0.4em -Q_{t_{2}} {S_{t_{2}}}^{-1}\hat{S}_{t_{2}} - (Q_{t_{2}} {S_{t_{2}}}^{-1}\hat{S}_{t_{2}} )^{'} + \hat{Q}_{t_{2}}
\end{align*}
where $(\hat{P},\hat{S},\hat{Q})_{t_{2}}$ are the rates of change of $(P,S,Q)_{t_{2}}$ under the 
dynamics at time $t$, parametrized by $(A(t), C(t), \Sigma(t))$ as per \er{eq:ric:biriccati}.
\begin{align*}
-\hat{P}_{t_{2}} &\doteq A(t){'}P_{t_{2}}+P A(t)+P_{t_{2}}\Sigma(t) P_{t_{2}} \\
-\hat{S}_{t_{2}} &\doteq (A(t)+\Sigma(t)P_{t_{2}}){'}S_{t_{2}} \\
-\hat{Q}_{t_{2}} &\doteq S_{t_{2}}{'}\Sigma(t)S_{t_{2}}
\end{align*}

This shows that the dual quadratic also satisfies a Riccati equation \er{eq:ric:ricd}
with coefficients
\begin{align}
\bar{A}(t)=&{S_{t_{2}}}^{-1}( \hat{P}_{t_{2}}{S_{t_{2}}^{-1}}^{'} Q_{t_{2}}- \hat{S}_{t_{2}}) \nonumber\\
\bar{\Sigma}(t)=&{S_{t_{2}}}^{-1}\hat{P}_{t_{2}}{S_{t_{2}}^{-1}}^{'} \label{eq:ric:dualriccoef}\\
\bar{C}(t)=& Q_{t_{2}} {S_{t_{2}}}^{-1}\hat{P}_{t_{2}}{{S_{t_{2}}}^{-1}}{'}Q_{t_{2}}  
 -Q_{t_{2}} {S_{t_{2}}}^{-1}\hat{S}_{t_{2}} - (Q_{t_{2}} {S_{t_{2}}}^{-1}\hat{S}_{t_{2}} )^{'} + \hat{Q}_{t_{2}} \nonumber
\end{align}
Above is equivalent to the following compatibility conditions which emphasize the symmetry and isolate $(\hat{P},\hat{S},\hat{Q})_{t_{2}}$.
\begin{align}\label{eq:ric:matchcond}
\begin{split}
-\hat{P}_{t_{2}}&=A(t){'}P_{t_{2}}+P_{t_{2}}A(t)+C(t)+P_{t_{2}}\Sigma(t) P=S_{t_{2}}\bar\Sigma(t) S_{t_{2}}{'}  \\
-\hat{S}_{t_{2}}&=(A(t)+\Sigma(t) P_{t_{2}}){'} S_{t_{2}}=S_{t_{2}}(-\bar{A}(t)+\bar\Sigma(t) Q_{t_{2}}) \\
-\hat{Q}_{t_{2}}&=S_{t_{2}}{'}\Sigma(t) S_{t_{2}}=-\bar{A}(t){'}Q_{t_{2}}-Q_{t_{2}}\bar{A}+\bar{C}(t)+Q_{t_{2}}\bar\Sigma(t) Q_{t_{2}}.
\end{split}
\end{align}
Through some algebraic manipulations, above compatibility conditions can be shown to
be equivalent to the matrix equality \er{eq:ric:hamilmatch}. Hence proved.
\end{proof}
\begin{remark}
Above result can be used to transform one DRE into its semiconvex dual, which maybe easier to solve.
For time invariant problems, it is possible to find the kernel to make $\bar{\Sigma}(t)\equiv 0$ and convert the dual DRE into a linear ODE, which permits analytical solutions. 
Algorithmically, first we compute the dual at the terminal time, 
$\mathcal{Q}_{t_{2}}=\mathcal{D}_{V_{t_{2}}}\mathcal{P}_{t_{2}}$. Then, the dual DRE is solved. Finally, solution to the original DRE can be computed using inverse dual operation,
$ \mathcal{P}_{t_{1}}=\mathcal{D}_{V_{t_{2}}}^{-1}\mathcal{Q}_{t_{1}}$.

Several known solutions to the time-invarient DRE e.g \cite{Rusnak}, \cite{Leipnik2}
can be derived from above general approach.
\end{remark}

\subsection{Kernel Relationships}

Next, we shall prove various kernel relationships between the primal and dual DRE solutions at different times.
This provides an insight into the solution manifolds and their structure. In process, we extend the \maxp solution described in
\cite{Funda} to the time-varying problem.

From \er{eq:ric:assm2}, we already have following kernel relationships. Namely
$\mathcal{D}_{V_{t_{2}}}[\mathcal{P}_{t_{2}}] = \mathcal{Q}_{t_{2}}$ and
$\mathcal{D}_{V_{t_{2}}}[\mathcal{P}_{t_{1}}] = \mathcal{Q}_{t_{1}}$.

\begin{theorem}\label{thm:kernelop1}
Assume \er{eq:ric:assm1} and \er{eq:ric:assm2}. Recall the bivariate quadratic \maxp kernel $I_{t_{1}}^{t_{2}}$ from \er{eq:ric:kernelform}. Then the following semiconvex duality relationships hold.
\begin{align}
\mathcal{D}_{I_{t_{1}}^{t_{2}}}[\mathcal{P}_{t_{1}}] = \mathcal{P}_{t_{2}} \label{eq:kernelop1:a}\\
\mathcal{D}_{V_{t_{1}}}[\mathcal{P}_{t_{1}}] = \mathcal{Q}_{t_{2}} \label{eq:kernelop1:b}
\end{align}
Thus we can propagate the DRE \er{eq:ric:introdre} from $\bar{p}_{t_{2}}$ to $\bar{p}_{t_{1}}$ using following kernel operations.
\begin{align}\label{eq:kernelop1:c}
\mathcal{P}_{t_{1}}=\mathcal{S}_{t_{1}}^{t_{2}}[\mathcal{P}_{t_{2}}]
=\mathcal{D}_{I_{t_{1}}^{t_{2}}}^{-1}[\mathcal{P}_{t_{2}}]
= \mathcal{D}_{V_{t_{1}}}^{-1} \mathcal{D}_{V_{t_{2}}} [\mathcal{P}_{t_{2}}]
\end{align}

\end{theorem}
\begin{proof}
Using the value function \er{eq:ric:simpleV} underlying the DRE \er{eq:ric:introdre}, teh semigroup operator \er{eq:ric:semigroup} and the dynamic programming principle, we get $\mathcal{P}_{t_{1}}(x)=\mathcal{S}_{t_{1}}^{t_{2}}[\mathcal{P}_{t_{2}}](x)$ for all $x\in\real^{n}$.
Since $I_{t_{1}}^{t_{2}}$ is the \maxp fundamental solution, using \er{eq:ric:kernelop} we have
$$ \mathcal{P}_{t_{1}}(x)=\mathcal{S}_{t_{1}}^{t_{2}}[\mathcal{P}_{t_{2}}]
= \sup_{y\in\real^{n}}\left( I_{t_{1}}^{t_{2}}(x,y)+ \mathcal{P}_{t_{2}}(y) \right)
= \mathcal{D}_{I_{t_{1}}^{t_{2}}}^{-1}[\mathcal{P}_{t_{2}}]. $$
This proves first half of \er{eq:kernelop1:c}.
Operating with $\mathcal{D}_{I_{t_{1}}^{t_{2}}}$ on both sides proves \er{eq:kernelop1:a}.

Using above result, \er{eq:ric:assm2} and \er{eq:ric:DPP} we have
\begin{align*}
\mathcal{P}_{t_{1}}(x) 
&=\mathcal{D}_{I_{t_{1}}^{t_{2}}}^{-1}[\mathcal{P}_{t_{1}}](x) = \mathcal{D}_{I_{t_{1}}^{t_{2}}}^{-1} \mathcal{D}_{V_{t_{2}}}^{-1}[\mathcal{Q}_{t_{2}}](x)\\
&=\sup_{y\in\real^{n}}\left( I_{t_{1}}^{t_{2}}(x,y)+\sup_{z\in\real^{n}}\left( V_{t_{2}}(y,z)+\mathcal{Q}_{t_{2}}(z) \right )\right)\\
&=\sup_{z\in\real^{n}}\left( \mathcal{Q}_{t_{2}}(z) +\sup_{y\in\real^{n}}\left( I_{t_{1}}^{t_{2}}(x,y)+V_{t_{2}}(y,z) \right )\right)\\
&=\sup_{z\in\real^{n}}\left( \mathcal{Q}_{t_{2}}(z) +\mathcal{S}_{t_{1}}^{t_{2}}[V_{t_{2}}(\cdot, z)](x) \right)\\
&=\sup_{z\in\real^{n}}\left( \mathcal{Q}_{t_{2}}(z) +V_{t_{1}}(x,z) \right)
=\mathcal{D}_{V_{t_{1}}}^{-1} [\mathcal{Q}_{t_{2}}](x)
=\mathcal{D}_{V_{t_{1}}}^{-1} [\mathcal{D}_{V_{t_{2}}}^{-1} \mathcal{P}_{t_{2}}](x)
\end{align*}
This proves later half of \er{eq:kernelop1:c}. 
Furthermore, operating by $\mathcal{D}_{V_{t_{1}}}$ on both sides proves \er{eq:kernelop1:b}.
Hence proved.
\end{proof}
\begin{remark}
Using \er{eq:ric:semicon} and \er{eq:ric:semiconinv}, it can be easily shown that above semiconvex dual relations lead to relations between the underlying parameters, some of which have been derived before. 
Specifically, equation \er{eq:kernelop1:a} leads to \er{eq:kernel_to_dresol:a},
equation \er{eq:kernelop1:b} leads to \er{eq:ric:maxpluspropag2}
which can also be derived from equation \er{eq:kernelop1:b} which implies
\be
\bar{q}_{t_{2}}=-S_{t_{2}}{'}(\bar{p}_{t_{2}}-P_{t_{2}})^{-1}S_{t_{2}}-Q_{t_{2}} =-S_{t_{1}}{'}(\bar{p}_{t_{1}}-P_{t_{1}})^{-1}S_{t_{1}}-Q_{t_{1}}. \label{eq:ric:propagp}
\ee

\end{remark}

Now we shall obtain a time-varying version of the result previously obtained in 
\cite{Funda}, in order to complete our picture of kernel relationships between primal and dual DREs.
%

\begin{theorem}\label{thm:ric:dualkernelop}
Assume \er{eq:ric:assm1} and \er{eq:ric:assm2}. $\mathcal{Q}_{t_{1}}(z)$ is the semiconvex primal of $\mathcal{Q}_{t_{2}}(\bar{z})$ under kernel $B_{t_{1}}^{t_{2}}(z,\bar{z}) $
\begin{align} 
\mathcal{Q}_{t_{1}}(z)&=\sup_{y\in\real^{n}}\left[ B_{t_{1}}^{t_{2}}(z,y) +\mathcal{Q}_{t_{2}}(y)\right]
=\mathcal{D}_{B_{t_{1}}^{t_{2}}}^{-1}[\mathcal{Q}_{t_{2}}](z) \label{eq:ric:dualkernelop1} \\
\intertext{where,} 
\label{eq:ric:dualkernelop2}
\mathcal{B}_{t_{1}}^{t_{2}}(z,y)&=\inf_{x\in\real^{n}}\left\{ V_{t_{1}}(x,y)- V_{t_{2}}(x,z) \right\}
\end{align}
Hence
$\mathcal{B}_{t_{1}}^{t_{2}}(z,y)=\frac{1}{2}z{'}{B_{11}}_{t_{1}}^{t_{2}}z+z{'}{B_{12}}_{t_{1}}^{t_{2}}y+\frac{1}{2}y{'}{B_{22}}_{t_{1}}^{t_{2}}y$, with
\begin{align}\label{eq:ric:dualkernelop3}
\begin{split}
{B_{11}}_{t_{1}}^{t_{2}}&=-S_{t_{2}}{'}(P_{t_{1}}-P_{t_{2}})^{-1}S_{t_{2}}-Q_{t_{2}}\\
{B_{12}}_{t_{1}}^{t_{2}}&=S_{t_{2}}{'}(P_{t_{1}}-P_{t_{2}})^{-1}S_{t_{1}}\\
{B_{22}}_{t_{1}}^{t_{2}}&=-S_{t_{1}}{'}(P_{t_{1}}-P_{t_{2}})^{-1}S_{t_{1}}-Q_{t_{1}}
\end{split}
\end{align}
and
\be \label{eq:ric:dualkernelop4}
\bar{q}_{t_{1}}={B_{11}}_{t_{1}}^{t_{2}}-{B_{12}}_{t_{1}}^{t_{2}}({B_{22}}_{t_{1}}^{t_{2}}+\bar{q}_{t_{2}})^{-1}{B_{12}}_{t_{1}}^{t_{2}}{'}
\ee
\end{theorem}
\begin{proof}
In theorem \ref{thm:kernelop1}, we proved that $\mathcal{P}_{t_{1}}=\mathcal{D}_{V_{t_{2}}}^{-1}[\mathcal{Q}_{t_{1}}]=\mathcal{D}_{V_{t_{1}}}^{-1}[\mathcal{Q}_{t_{2}}]$. 
Hence we have
\begin{align}
\mathcal{Q}_{t_{1}}(z) &= \mathcal{D}_{V_{t_{2}}}\mathcal{D}_{V_{t_{1}}}^{-1}[\mathcal{Q}_{t_{2}}](z)\nonumber\\
&= \inf_{x\in\real^{n}}\left( \mathcal{D}_{V_{t_{1}}}^{-1}[\mathcal{Q}_{t_{2}}](x) - V_{t_{2}}(x,z) \right) \nonumber\\
&= \inf_{x\in\real^{n}}\left( \sup_{y\in\real^{n}}\left(V_{t_{1}}(x,y)+\mathcal{Q}_{t_{2}}(y)\right) - V_{t_{2}}(x,z) \right) \nonumber\\
&= \inf_{x\in\real^{n}} \sup_{y\in\real^{n}}\left(V_{t_{1}}(x,y)+\mathcal{Q}_{t_{2}}(y)- V_{t_{2}}(x,z)\right) \nonumber\\
&= \inf_{x\in\real^{n}} \sup_{y\in\real^{n}}\psi(x,y) \label{eq:ric:dualkernelop:a}
\end{align}
where 
$$\psi(x,y)=\frac{1}{2}x{'}(P_{t_{1}}-P_{t_{2}})x+x{'}S_{t_{1}}y+\frac{1}{2}y{'}(Q_{t_{2}}+\bar{q}_{t_{2}})y-x{'}S_{t_{2}}z-\frac{1}{2}z{'}Q_{t_{2}}z$$
Note that by \er{eq:ric:assm2}, $P_{t_{1}}-P_{t_{2}}\succ 0$. Hence $\psi(x,y)$ is strictly convex in $x$. Also observe that by theorem \ref{thm:kernelop1}, $\mathcal{P}_{t_{1}}(x)=\sup_{y}(V_{t_{1}}(x,y)+\mathcal{Q}_{t_{2}}(y))$ exists  for any $x\in \real^{n}$. Hence $Q_{t_{2}}+\bar{q}_{t_{2}}\prec 0$. Thus $\psi(x,y)$ is strictly concave in $y$. For such a convex-concave function following saddle point exists. By setting $\grad_{x}\psi$ and $\grad_{y}\psi$ equal to zero, and solving, we get 
\begin{align*}
x^{0}&=\left((P_{t_{1}}-P_{t_{2}})-S_{t_{1}}(Q_{t_{2}}+\bar{q}_{t_{2}})^{-1}S_{t_{1}}{'}\right)^{-1}S_{t_{2}}z\\
y^{0}&=-(Q_{t_{2}}+\bar{q}_{t_{2}})^{-1}S_{t_{1}}{'}x^{0}
\end{align*}
For such $x^{0}$ and $y^{0}$, $\psi(x^{0},y)\leq \psi(x^{0},y^{0})\leq \psi(x,y^{0})$. 
Hence by a well known result,
\be \label{eq:ric:dualkernelop:b}
\inf_{x\in\real^{n}} \sup_{y\in\real^{n}}\psi(x,y)=\psi(x^{0},y^{0})= \sup_{y\in\real^{n}}\inf_{x\in\real^{n}}\psi(x,y)
\ee
Using \er{eq:ric:dualkernelop:a} and \er{eq:ric:dualkernelop:b}
\begin{align}
\mathcal{Q}_{t_{1}}(z)&=\sup_{y\in\real^{n}}\inf_{x\in\real^{n}}\psi(x,y)\\
&=\sup_{y\in\real^{n}}\left( \mathcal{Q}_{t_{2}}(y) +\inf_{x\in \real^{n}}(V_{t_{1}}(x,y)-V_{t_{2}}(x,z))\right)\\
&=\sup_{y\in\real^{n}}\left( \mathcal{Q}_{t_{2}}(y) +B_{t_{1}}^{t_{2}}(z,y)\right)
\end{align}
\er{eq:ric:dualkernelop3} can be easily obtained from \er{eq:ric:dualkernelop2} by finding local minimum in $x$ (which is global minimum, since infimum exists), substituting and term-wise equating coefficients. Similarly \er{eq:ric:dualkernelop4} results from substituting $\mathcal{Q}_{t}=\frac{1}{2}z{'}\bar{q}_{t}z$ , \er{eq:ric:dualkernelop2}, \er{eq:ric:dualkernelop3} into \er{eq:ric:dualkernelop1}.
\end{proof}

Thus the asssumptions \er{eq:ric:assm2} and theorems \ref{thm:kernelop1} and \ref{thm:ric:dualkernelop} 
can be summarized in the diagram below. Note that primal and dual quadratics are on top and bottom respectively. Vertical and diagonal lines show duality transformation with indicated kernel. Arrows are directed from the primal to it semiconvex dual.
\begin{figure}[h]
\label{figb} 
\centering
\begin{picture}(180,120)
\thicklines
\put(-20,120){Primal DRE: $-\dot{\bar{p}}_{t}=A(t){'}\bar{p}_{t}+\bar{p}_{t}A(t)+C(t)+\bar{p}_{t}\Sigma(t)\bar{p}_{t}$.}
\put(10,90){$\mathcal{P}_{t_2}(x)$} 
\put(160,90){$\mathcal{P}_{t_1}(x)$} 
\put(10,20){$\mathcal{Q}_{t_2}(z)$} 
\put(160,20){$\mathcal{Q}_{t_1}(z)$}

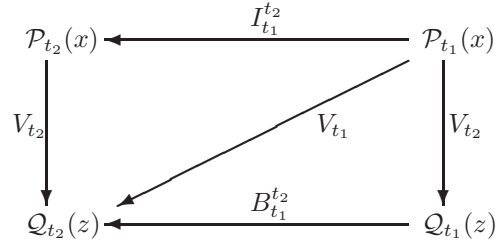
\put(-20,-10){Dual DRE: $-\dot{\bar{q}}_{t}=\bar{A}(t){'}\bar{q}_{t}+\bar{q}_{t}\bar{A}(t)+\bar{C}(t)+\bar{q}_{t}\bar{\Sigma}(t)\bar{q}_{t}$.}

\put(155,93){\vector(-1,0){115}}

\put(18,85){\vector(0,-1){55}}

\put(155,23){\vector(-1,0){115}}

\put(168,85){\vector(0,-1){55}}

\put(155,85){\vector(-2,-1){110}}

\put(5,58){$V_{t_2}$}
\put(170,58){$V_{t_2}$}
\put(120,58){$V_{t_1}$}

\put(95,98){$I_{t_1}^{t_2}$}
\put(95,28){$B_{t_1}^{t_2}$}
\end{picture}
\vskip 1em
\caption{Time varying problem: Duality relationships.}
\end{figure}

In summary, so far we saw three distinct ways of solving \er{eq:ric:introdre}, that is obtaining $p_{t_{1}}$ from $p_{t_{2}}$ all of which can be expressed using the semiconvex duality.
\begin{enumerate}
\item Direct method which  assumes only \er{eq:ric:assm1}. Formulae are given by  \er{eq:kernel_to_dresol:a} and \er{eq:ric:kernelform}. Propagation is achieved by following transform.
$$\mathcal{P}_{t_{1}}=\mathcal{D}_{I_{t_{1}}^{t_{2}}}^{-1}[\mathcal{P}_{t_{2}}]$$
Problem with this method is that as $t_{1}\rightarrow t_{2}$, parameters of the kernel $I_{t_{1}}^{t_{2}}$ become singular, limiting solution accuracy. 
\item Alternate method, which assumes \er{eq:ric:assm1} and \er{eq:ric:assm2}. Formulas are given by \er{eq:ric:propagp}, which is same as \er{eq:ric:maxpluspropag2}. Propagation is achieved by following transform. 
$$\mathcal{P}_{t_{1}}=\mathcal{D}_{V_{t_{1}}}^{-1}\mathcal{D}_{V_{t_{2}}}[\mathcal{P}_{t_{2}}]$$
This method works better for a small time step propagation, since parameters of kernels $V_{t_{1}}$ and $V_{t_{2}}$ do not blow up.
\item Third method also assumes \er{eq:ric:assm1} and \er{eq:ric:assm2}. Time invariant version of this method was first proposed in \cite{Funda}. Formulae are \er{eq:ric:dualkernelop3} and \er{eq:ric:dualkernelop4}. Propagation is achieved by following transform
$$\mathcal{P}_{t_{1}}=\mathcal{D}_{V_{t_{2}}}^{-1}\mathcal{D}_{B_{t_{1}}^{t_{2}}}^{-1}\mathcal{D}_{V_{t_{2}}}[\mathcal{P}_{t_{2}}]$$

Problem with this method is similar to the direct method. Namely, $t_{1}\rightarrow t_{2}$, parameters of the kernel $I_{t_{1}}^{t_{2}}$ become more singular, limiting solution accuracy.
\end{enumerate}

\section{Conclusion}
In conclusion, this paper derives a new representation of the solution to the time-varying DRE \er{eq:ric:introdre} based on the \maxp fundamental solution and semiconvex duality. It derives the \maxp kernel and the fundamental solution first
proposed in \cite{FlemMac} for the time varying linear quadratic problem. Such solution also solves a two point optimal control problem \er{eq:ric:Fundakernel}.
Equations  \er{eq:ric:kernelform} describe the kernel, which is bivariate quadratic in terms of the two boundary points, and its parameters can computed using the evolution of three matrix ODEs \er{eq:ric:biriccati} or the transition matrix of the associated Hamiltonian matrix \er{eq:ric:IfromPhi}. The formulae for the time evolution and combination of the kernel are also derived in 
\er{eq:ric:biriccatiI} and \er{eq:ric:combop2}. Fundamental solution can be used to solve the DRE analytically, starting from any initial condition, using \er{eq:kernel_to_dresol:a} and \er{eq:ric:maxpluspropag2}.

This paper also shows that under a fixed duality kernel, semiconvex dual of the DRE solution satisfies another dual DRE, which coefficient satisfy the matrix compatibility conditions \er{eq:ric:hamilmatch} involving Hamiltonian and certain symplectic matrices.
For time invariant DRE, this allows us to make dual DRE linear and thereby solve the primal DRE analytically. 

This paper also derives various kernel/duality relationships between the primal and dual DREs at different times shown in figure \ref{figb},
which leads to an array of methods to solve a DRE. Time invariant analogue of one of these methods was first proposed in \cite{Funda}.

Thus this paper provides a \maxp fundamental solution for the time-varying linear quadratic DRE, useful for stiff problems and long time horizon propagation. It also provides a powerful unifying framework based on optimal control formulation, semiconvex duality and \maxp algebra, which enables us to solve the Riccati differential equations,
and see existing methods in new light.

 \bibliographystyle{halpha}  
 \bibliography{myrefs}  

\appendix
\section{Appendix}
\begin{lemma}
\label{lemma:biquadhjbexist1}
Define 
\be \label{eq:ric:Wzdef}
W(t,x;z)\doteq \frac{1}{2}x{'}P_t x+x{'}S_t z+\frac{1}{2}z{'}Q_t z,
\ee
where $P_{t}$, $Q_{t}$, $S_{t}$ evolve according to the \er{eq:ric:biriccati}. 
Then $S_{t}$ is invertible for $t\in (\bar{T}, T]$
and $W(t,x;z)$
is the solution of the following Hamilton-Jacobi-Bellman PDE on $\left(\bar{T},T\right]\times \real^n$.
\begin{align}
 0&=-\grad_{t}{W}(t,x; z)-H\left(t,x,\grad_{x} W(t,x;z)\right), \label{eq:ric:hjbpde} \\
 \intertext{where the Hamiltonian $H$ is defined as }
H(t,x,p)& \doteq \sup_{u\in\real^n}\left\{ p{'}f_t(x,u)+l_t(x,u) \right\} \nonumber\\
 &= \sup_{u\in\real^n}\left\{ p{'}\left(A(t)x+\sigma(t)u\right)+\frac{1}{2}x{'}C(t)x-\frac{1}{2}|u|^2\, dt \right\} \nonumber\\
&=\frac{1}{2}x{'}C(t)x+x{'}A{'}(t)p+\frac{1}{2}p{'}\Sigma(t) p, 
\label{eq:ric:hjbpde2} 
\end{align}
and with the terminal payoff defined in \er{eq:biric:payoffeq} as the boundary condition at time $T$,
\be
W(T,x;z)=\phi(x,z)\doteq \frac{1}{2}x{'} P_{T} x+x{'}S_{T} z +\frac{1}{2}z{'}Q_{T} z \hskip 2em \forall x \in \real^n. \label{eq:ric:bddcond}
\ee

\end{lemma}
\begin{proof}
Existence of the solution $P_{t}:-\bar{T}<t\leq T$ is assumed in \er{eq:ric:assm1}.
This combined with local boundedness, and  piecewise continuity of coefficients
guarantees the existence of $S_{t}$, and hence that of $Q_{t}$ for 
$-\bar{T}<t\leq T$.

Now we prove that $S_{t}$ is invertible.
Let us define, $B(t)=-(A(t)+\Sigma(t)P_t)$. Then $S_{t_1}=\Phi_B(t_1,T)S_{T}$, 
where $\Phi_B$ is the state transition matrix of the system $\dot{x}_t=B(t)x_t$. By Abel-Jacobi-Liouville formula
$$\det \Phi_B(t_1,T)=e^{\int_{t_1}^{T}\mathbf{Tr} B(\tau)\,d\tau} > 0$$
Since both $\Phi_B(t_1,T)$ and $S_T=S$ are invertible,
$S_{t_1}=\Phi_B(t_1,T)S_{T}$ is invertible as well.

The proof that $W$ solves HJB PDE, is immediate by substitution in \er{eq:ric:hjbpde} and \er{eq:ric:bddcond}.

\end{proof}

Next we need a verification theorem to connect HJB PDE solution to the control value function.
\begin{lemma}
\label{lem:ric:verification}
Assume \er{eq:ric:assm1}. Let $W$ and $J$ be defined as per \er{eq:ric:Wzdef} and \er{eq:biric:payoffeq}, respectively. Let $x,z\in \real^n$ and $t_1 \in \left(\bar{T},T\right]$. One has
$$W(t,x;z)\geq {J_{t}^{T}}(x,u;z) \hskip 2em \forall u\in L_2(t,T)$$
and
$$W(t,x;z)={J_{t}^{T}}(x,\tilde{u};z)$$
where $\tilde{u}_{s}=\tilde{u}(s,\xi_{s})=\sigma(s){'}\grad W(s,\xi_{s};z)
=\sigma(s){'}(P_{s} \xi_{s}+S_{s} z)$, which implies $W^z=V^z$ and
\be
V_{t}(x;z)=W_{t}(x;z)=\frac{1}{2}x{'}P_t x+x{'}S_t z+\frac{1}{2}z{'}Q_t z
\label{eq:ric:valuesolveshjb}\ee
\end{lemma}
\begin{proof}
Let $u\in L_2(t,T)$.
\begin{align*}
{J_{t}^{T}}(x,u;z){}=&\int_{t}^{T} \left(L_{s}(x_{s},u_{s})+(A(s)x_{s}+\sigma(s)u_{s}){'}\grad W(s,\xi_{s};z)\right)\, ds +\phi(\xi_{T};z)\\
{}& -\int_{t}^{T} (A(s)\xi_{s}+\sigma(s)u_{s}){'}\grad W(s,\xi_{s};z)\, ds\\
\intertext{which by definition of $H$}
{}\leq&\int_{t}^{T} H\left(\xi_t,\grad W(s,\xi_{s};z)\right)\, ds+\phi(\xi_{T},z)-
\int_{t}^{T} (A(s)\xi_{s}+\sigma(s)u_{s}){'}\grad W(s,\xi_{s};z)\, ds \\
\intertext{which by \er{eq:ric:hjbpde} and \er{eq:ric:dyneq}}
{}=&\int_{t}^{T}\left\{ -\frac{\partial}{\partial s} W(s,\xi_{s};z)-\dot{\xi}_s\grad W(s,\xi_{s};z)\right\}\, ds+ \phi(\xi_{T};z)\\
{}=& -\int_{t}^{T}\frac{d}{dt}W(s,\xi_{s};z)\, ds +\phi(\xi_{T};z)\\
{}=&W(t,x;z)-W(T,\xi_{T};z)+\phi(\xi_{T};z)=W(t,x;z)
\end{align*}
using \er{eq:ric:bddcond}.

Also note that in the above proof, if we substitute $\tilde{u}_{s}=\sigma(s){'}\grad W(s,\xi_{s};z)$, then using $l_{s}(x,u)=\frac{1}{2}x{'}C(s)x - |u_{s}|^{2}/2$,
$\sigma(s)\sigma{'}(s)=\Sigma(s)$ and \er{eq:ric:hjbpde2},
\begin{align*}
&l_{s}(\xi_{s},u_{s})+(A(s)\xi_{s}+\sigma(s)u_{s}){'}\grad W(s,\xi_{s};z)\\
& {}= \frac{1}{2}\xi_{s}{'}C(s)\xi_{s}+\frac{1}{2}\grad {W}{'}(s,\xi_{s};z)\Sigma(s)\grad {W}(s,\xi_{s};z)+A{'}(s)\grad {W}(s,\xi_{s};z)\\
&{}=H\left(s,\xi_{s},\grad {W}(s,\xi_{s};z)\right) 
\end{align*}
This converts the inequality into equality and we get ${J_{t}^{T}}(x,\tilde{u};z)=W(t,x;z)$. 
\end{proof}


First we derive a lemma about the end point of optimal trajectories.

\begin{lemma} \label{lem:ric:opttrajend}
Assume \er{eq:ric:assm1}. Consider the system trajectory $\tilde{\xi}$ evolving according to \er{eq:ric:dyneq} under optimal control $\tilde{u}_{s}=\sigma(s){'}(P_{s}\tilde{\xi}_{s}+S_{s}z)$ from lemma \ref{lem:ric:verification}. Then for $\bar{T}<t_{1}\leq t_{2}\leq T$,

$$ 
S_{t_{2}}{'}\tilde{\xi}_{t_{2}}+Q_{t_{2}}z=S_{t_{1}}{'}\tilde{\xi}_{t_{1}}+Q_{t_{1}}z
$$

\end{lemma}
\begin{proof}
By time-varying linear system theory, for a system evolving as per 
\begin{align*}
\dot{\tilde{\xi}}_{s}&=A(s)\tilde{\xi}_{s}+\sigma(s)\tilde{u}_{s}\\
&{}=A(s)\tilde{\xi}_{s}+\sigma(s)\sigma(s){'}(P_{s}\xi_{s}+S_{s}z) \\
&{}=\left( A(s) + \Sigma(s)P_{s}\right) \xi_{s}+\Sigma(s)S_{s}z
\end{align*}
solution is given as
\be
\tilde{\xi}_{t_{2}}=\Phi_{B}(t_{2},t_{1})\tilde{\xi}_{t_{1}}+\int_{t_{1}}^{t_{2}}\Phi_{B}(t_{2},s)\Sigma(s)S_{s}z\, d\tau \label{eq:ric:opttrajend:a}
\ee
where $\Phi_{B}(t_{2},t_{1})=U_{t_{2}}U_{t_{1}}^{-1}$ and $U$ is the solution of differential equation $\dot{U_{s}}=B(s)U_{s}$, with $B(s)=A(s)+\Sigma(s)P_{s}$. 

It is well known that the state transition matrix 
$$\Phi_{B(s)}(t_{2},t_{1})=\Phi_{-B(s)'}'(t_{1},t_{2})$$
now, noting from \er{eq:ric:biriccati} that $\dot{S}_{s}=-(A(s)+\Sigma(s)P_{s}){'}S_{s}=-B(s){'}S_{s}$, and since $S_{t_2}$ is invertible, we have
\be\label{eq:ric:opttrajend:b}
\Phi_{B(s)}(t_{2},t_{1})=\Phi_{-B(s)'}'(t_{1},t_{2})=\left(S_{t_{1}}S_{t_{2}}^{-1} \right){'}={S_{t_{2}}^{-1}}{'}S_{t_{1}}{'}
\ee
Substituting in \er{eq:ric:opttrajend:a}, and noting from \er{eq:ric:biriccati} that $-\dot{Q}_{s}=S_{s}{'}\Sigma(s)S_{s}$, 
\begin{align*}
\tilde{\xi}_{t_{2}}&={S_{t_{2}}^{-1}}{'}S_{t_{1}}{'}\tilde{\xi}_{t_{1}}+{S_{t_{2}}^{-1}}{'}\int_{t_{1}}^{t_{2}}S_{\tau}{'}\Sigma(\tau)S_{\tau}z\, d\tau \\
&={S_{t_{2}}^{-1}}{'}S_{t_{1}}{'}\tilde{\xi}_{t_{1}}+{S_{t_{2}}^{-1}}{'}\left(\int_{t_{1}}^{t_{2}}S_{\tau}{'}\Sigma(\tau)S_{\tau}\, d\tau \right)z \\
&={S_{t_{2}}^{-1}}{'}S_{t_{1}}{'}\tilde{\xi}_{t_{1}}+{S_{t_{2}}^{-1}}{'}\left(
Q_{t_{1}}-Q_{t_{2}}\right)z
\end{align*}
thus we have,
$$S_{t_{2}}{'}\tilde{\xi}_{t_{2}}+Q_{t_{2}}z=S_{t_{1}}{'}\tilde{\xi}_{t_{1}}+Q_{t_{1}}z $$
\end{proof}

Now we shall prove another useful lemma before turning to the main result.
\begin{lemma}\label{lem:ric:posdefQ}
Given $\bar{T}<{t_1}<{t_2}\leq T$, and $Q_t$ evolving according to \er{eq:ric:biriccati} with terminal value $Q_T=Q$, then
$$ Q_{t_1}-Q_{t_2} \succ 0$$
\end{lemma}
\begin{proof}
Note that we assumed in \er{eq:ric:assm1} that the system $\dot{\xi}_{s}=A(s)\xi_{s}+\sigma(s)u_{s}$ parametrized by $(A(s),\sigma(s))$ is controllable. This is true if and only if the following {\it controllability grammian} is invertible for any $\bar{T}<{t_1}<{t_2}\leq T$.
\be\label{eq:ric:controlgram}
\int_{t_{1}}^{t_{2}}\Phi_{A}(t_{1},s)\sigma(s)\sigma(s){'}\Phi_{A}(t_{1},s){'}\, ds \succ 0
\ee

Thus for all $\xi,y\in \real^{n}$, $\exists$ control $\hat{u}_{s}$ such that is the trajectory $\dot{\hat{\xi}}=A(s)\hat{\xi}_{s}+\sigma(s)\hat{u}_{s}$ with $\hat{\xi}_{t_{1}} =x$ satisfies $\hat{\xi}_{t_{2}}=y$.

Now we claim that system $\left(A(s)+\Sigma(s)P_{s}, \sigma(s)\right)$ is also controllable. This is clear because by using control $\bar{u}_{s}= \hat{u}_{s}-\sigma(s){'}P_{s}\xi_{s}$, we can keep the system trajectory same and reach from $x$ to $y$.
\begin{align*}
\dot{\hat{\xi}}&=A(s)\hat{\xi}_{s}+\sigma(s)\hat{u}_{s}\\
&=(A(s)+\sigma(s)\sigma(s){'}P_{s})\xi_{s}+\sigma(s)\left( \hat{u}_{s}-\sigma(s){'}P_{s}\xi_{s}\right)\\
&=(A(s)+\Sigma(s)P_{s})\xi_{s}+\sigma(s)\bar{u}_{s}
\end{align*}
Hence similar to \er{eq:ric:controlgram}, using $B(s)=A(s)+\Sigma(s)P_{s}$ and $\sigma(s)\sigma(s){'}=\Sigma(s)$, the following controllability grammian is invertible.
\be\label{eq:ric:controlgram2}
\int_{t_{1}}^{t_{2}}\Phi_{B}(t_{1},s)\Sigma(s)\Phi_{B}(t_{1},s){'}\, ds \succ 0
\ee
Substituting $\Phi_{B}(t_{1},s)=S_{t_{1}}^{-1}{'}S_{s}{'}$ from \er{eq:ric:opttrajend:b},
\begin{align}
\int_{t_{1}}^{t_{2}}\Phi_{B}(t_{2},s)\Sigma(s)\Phi_{B}(t_{2},s){'}\, ds 
&=\int_{t_{1}}^{t_{2}}S_{t_{1}}^{-1}{'}S_{s}{'}\Sigma(t)S_{s}S_{t_{1}}^{-1}\, ds\nonumber\\
&=S_{t_{1}}^{-1}{'}\left\{\int_{t_{1}}^{t_{2}}S_{s}{'}\Sigma(t)S_{s}\, ds\right\}S_{t_{1}}^{-1}\nonumber\\
&=S_{t_{1}}^{-1}{'}\left(Q_{t_{1}}-Q_{t_{2}}\right)S_{t_{1}}^{-1} \label{eq:ric:controlgram3}
\end{align}
where in the last equation, we used $Q_{t}$ evolution from \er{eq:ric:biriccati}. Using \er{eq:ric:controlgram2} and since $S_{t_{1}}$ is invertible by Lemma \er{lemma:biquadhjbexist1},  we have $Q_{t_{1}}-Q_{t_{2}}\succ 0$. 
\end{proof}

\begin{theorem}\label{proof:thm:ric:semicon}
Below is the proof of theorem \ref{thm:ric:semicon}
\end{theorem}
\begin{proof}
Since $\mathcal{P} \in {\mathcal{S}}^{-P}$, $\mathcal{P}(x)-\phi(x,z)$ is convex in $x$. Now,
\begin{align*}
\sup_z[\phi(x,z)+\mathcal{Q}(z)]&= \sup_z[\phi(x,z)+\inf_{y}[\mathcal{P}(y)-\phi(y,z)]\\
&= \sup_z \inf_y [\mathcal{P}(y)+\phi(x,z)-\phi(y,z)]\\
&= \sup_z \inf_y [\mathcal{P}(y)+\frac{1}{2}x{'} Px-\frac{1}{2}y{'} Py+(x-y){'}Sz]\\
\intertext{Let $\bar{z}=Sz$. Since $S$ is invertible, $\bar{z}$ also spans $\real^n$}
&= \frac{1}{2}x{'} Px+ \sup_{\bar{z}}\inf_y[\mathcal{P}(y)-\frac{1}{2}y{'} Py+(x-y){'}\bar{z}]\\
&= \frac{1}{2}x{'} Px +\sup_{z'}[x{'}\bar{z}
+\inf_y[\mathcal{P}(y)-\frac{1}{2}y{'} Py -y{'}\bar{z}]\\
\intertext{Since $\mathcal{P}(y)-\frac{1}{2}y{'} Py$ is convex, by Legendre-Fenchel transform (see e.g. Theorem $11.1$ in \cite{rockafellar})}
&= \frac{1}{2}x{'} Px+\mathcal{P}(x)-\frac{1}{2}x{'} P x\\
&= \mathcal{P}(x)
\end{align*}
\end{proof}

\end{document}